\documentclass[12pt,reqno]{amsart}
\usepackage{amsmath}
\usepackage{amssymb}
\usepackage{verbatim}
\usepackage{mathrsfs}
\usepackage{graphicx}
\usepackage{caption}
\usepackage{subcaption}
\usepackage{epstopdf}
\usepackage{makecell}
\usepackage{tabularx}
\usepackage{longtable}
\usepackage{tikz}
\usepackage{pgfplots}
\usepackage{xcolor}
\usepackage{pgfplots}
\usepackage{pgfplotstable}
\usepackage[ruled,vlined]{algorithm2e}
\usepackage[noend]{algpseudocode}
\usepackage[toc]{appendix}
\usepackage[capitalise]{cleveref}
\usepackage[noadjust]{cite}
\usepackage{color}

\usetikzlibrary{arrows}

\usepackage[top=0.6in,bottom=0.5in,left=0.7in,right=0.7in]{geometry}

\newtheorem{lemma}{Lemma}

\newtheorem{theorem}[lemma]{Theorem}
\newtheorem{remark}[lemma]{Remark}

\newtheorem{example}{Example}

\numberwithin{equation}{section}
\numberwithin{lemma}{section}

\newcommand{\bo}{\mathcal{O}}

\newcommand{\ka}{\textsf{k}}

\newcommand{\be}{ \begin{equation} }
\newcommand{\ee}{ \end{equation} }

\begin{document}
		
\title[Distinct Numerical Solutions for Elliptic Cross-Interface Problems  Using FEM and FDM]{
Distinct Numerical Solutions for Elliptic Cross-Interface Problems Using Finite Element and Finite Difference Methods}

\author{Qiwei Feng}

	\thanks{This research was partially supported by the Mathematics Research Center, Department of Mathematics, University of Pittsburgh, Pittsburgh, PA, USA}

\address{Department of Mathematics, University of Pittsburgh, Pittsburgh, PA 15260, USA.
	\quad {\tt qif21@pitt.edu}}

\address{Department of Mathematical and Statistical Sciences,  University of Alberta, Edmonton, Alberta, T6G 2G1, Canada.  	\quad {\tt qfeng@ualberta.ca} }

\begin{abstract}

In this paper, we discuss the second-order finite element method (FEM) and finite difference method (FDM) for numerically solving elliptic cross-interface problems characterized by vertical and horizontal straight lines, piecewise constant coefficients, two homogeneous jump conditions, continuous source terms, and Dirichlet boundary conditions. For brevity, we consider a 2D simplified version where the intersection points of the interface lines coincide with grid points in uniform Cartesian grids. Our findings reveal interesting and important results:
(1) When the coefficient functions exhibit either high jumps with low-frequency oscillations or low jumps with high-frequency oscillations, the finite element method and finite difference method yield similar numerical solutions.
(2) However, when the interface problems involve high-contrast and high-frequency coefficient functions, the numerical solutions obtained from the finite element and finite difference methods differ significantly.
Given that the widely studied SPE10 benchmark problem (see https://www.spe.org/web/csp/datasets/set02.htm) typically involves high-contrast and high-frequency permeability due to varying geological layers in porous media, this phenomenon warrants attention. To our best knowledge, so far there is no available literature that has clearly observed such significant differences in the numerical solutions produced by finite element and finite difference methods.
Furthermore,
this observation is particularly important for developing multiscale methods, as reference solutions for these methods are usually obtained using the standard second-order finite element method with a fine mesh, and analytical solutions are not available.
We provide sufficient details to enable replication of our numerical results, and the implementation is straightforward. This simplicity ensures that readers can easily confirm the validity of our findings.

\end{abstract}


\maketitle

\maketitle

\pagenumbering{arabic}

\section{Introduction}\label{introdu:1:intersect}
\begin{figure}[htbp]
	\centering	
	\begin{subfigure}[b]{0.45\textwidth}
	\includegraphics[width=7.5cm,height=7.5cm]{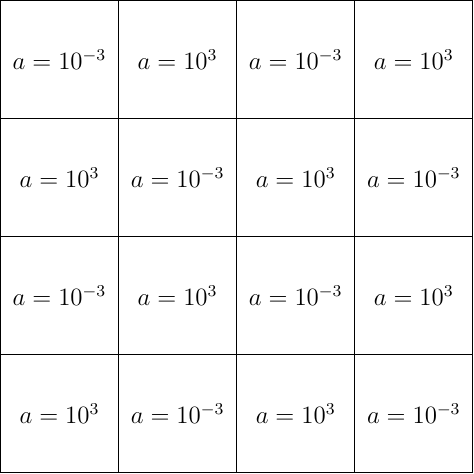}
\end{subfigure}
\hspace{0.5cm}
\begin{subfigure}[b]{0.45\textwidth}
	\includegraphics[width=7.5cm,height=7.5cm]{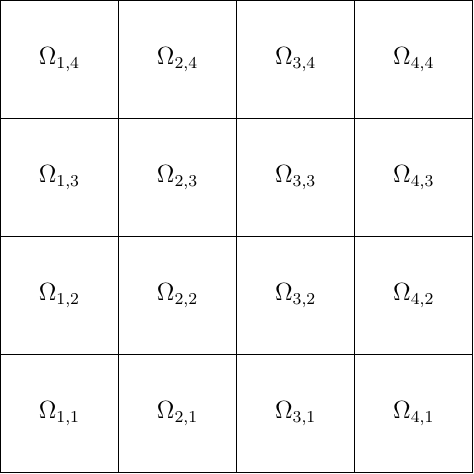}
\end{subfigure}
	\caption{An example for the coefficient function $a$ (left) and subdomains $\Omega_{p,q}$ (right) for the model problem in \eqref{intersect:3} with $m=4$.}
	\label{coefficient}
\end{figure}

Intersecting interface problems arise in many applications, such as the simulation of fluid flow in heterogeneous porous media  (see e.g. \cite{Minev2018,AliMankad2020,GuiraAusas2019,ArbogastXiao2013,ArbogastTao2013,Jaramillo2022,Kippe2008,Tahmasebi2018,Butler2020,TArbHXiao2015,Rasaei2008}). A well-known example of such a problem is the SPE10 benchmark problem, developed by the Society of Petroleum Engineers (see https://www.spe.org/web/csp/datasets/set02.htm). 
According to \cite{Vazqu07}, the high-frequency permeability in the various geological layers  of the SPE10 problem results in a highly oscillatory and high-contrast coefficient function of the interface problem. The singularity of the solution of the intersecting interface problem  is notably intensified when the permeability coefficient exhibits substantial discontinuities across interfaces, particularly when these discontinuities span multiple orders of magnitude (see e.g. \cite{CaKi01,Kell75,KiKo06,Blumen85,Petz01,Petz02,Nic90,KCPK07,NS94,FGW73,Kell71,Kell72}).
In this paper, we consider a simplified 2D variant of the elliptic intersecting interface problem, where the interfaces intersect along vertical and horizontal straight lines (see \cref{coefficient} for an illustration). Since an analytical solution typically does not exist even for this simplified problem, the reference solution is usually obtained by using standard second-order finite element method with linear basis functions and an appropriately fine mesh, such as various multiscale  methods (see e.g. \cite{Hellmanvist2017,PietroErn2012,lqvistPeterseim2014,EngwerHenning2019,LiHu2021,FuChungLi2019,HouHwang2017,ZhangCiHou2015,lqvistPersson2018}). However, our numerical results show that FEM and FDM produce markedly different solutions for elliptic intersecting interface problems with high-contrast and highly oscillatory coefficient functions (see \cref{fig:exam3,fig:exam4,fig:exam5,fig:exam7} in \cref{Intersect:ex3,Intersect:ex4,Intersect:ex5,Intersect:ex7}).

On the other hand, several approaches have been developed to address interface problems with smooth non-intersecting interfaces, piecewise smooth coefficients, two non-homogeneous jump conditions, discontinuous source terms, and mixed boundary conditions. Notable among these are immersed interface methods (IIM, see e.g. \cite{CFL19,LeLi94,EwingLLL99,GongLiLi08,HeLL2011,Li98,WieBube00,PanHeLi21,XiaolinZhong07,DFL20}), matched interface and boundary (MIB) methods (see e.g. \cite{YuZhouWei07,ZW06,FendZhao20pp109677,YuWei073D,ZZFW06}), and sixth-order schemes  (see e.g. \cite{FengHanMinevPOISS,FHM21Helmholtz,FengHanMinevHYB}).

%
%

\begin{figure}[htbp]
	\centering	
	\begin{subfigure}[b]{0.45\textwidth}
	\includegraphics[width=7.5cm,height=7.5cm]{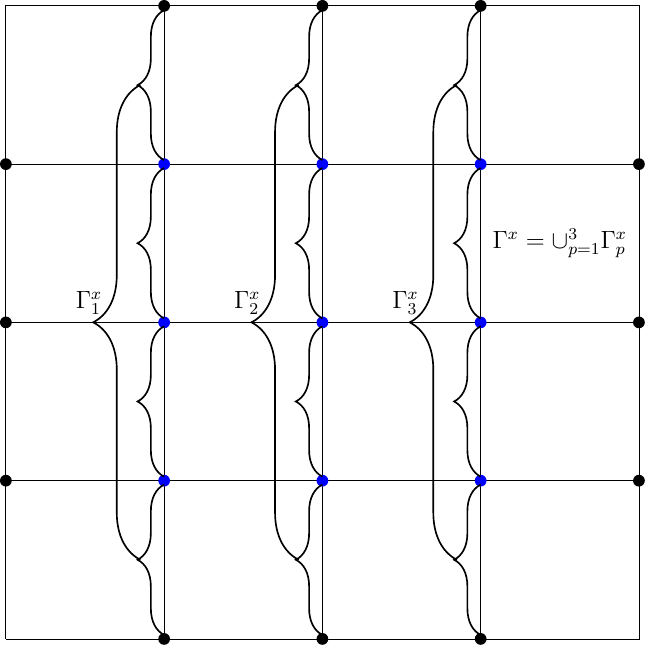}
\end{subfigure}
\hspace{0.5cm}
\begin{subfigure}[b]{0.45\textwidth}
	\includegraphics[width=7.5cm,height=7.5cm]{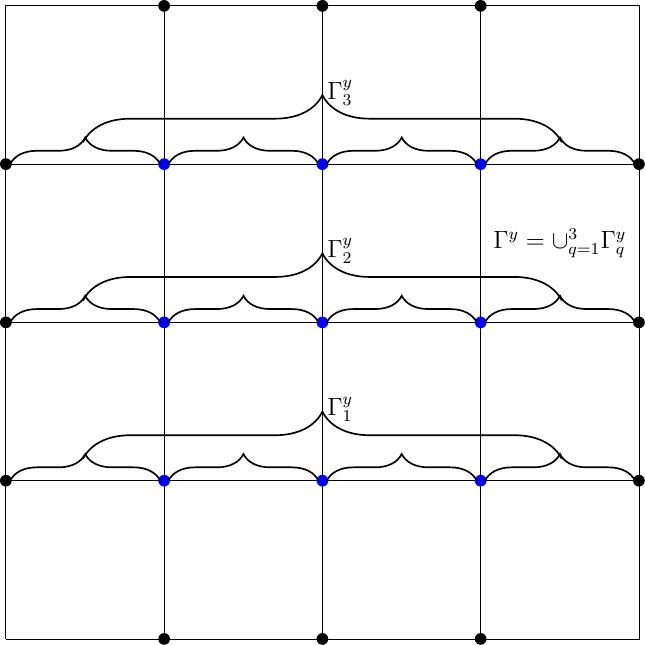}
\end{subfigure}
	\caption{An illustration for the vertical lines $\Gamma^x$ in \eqref{GammaXY} of the cross-interface $\Gamma$ for the model problem  \eqref{intersect:3} with $m=4$ (left). An illustration for the horizontal lines $\Gamma^y$ in \eqref{GammaXY} of the cross-interface $\Gamma$ for the model problem \eqref{intersect:3} with $m=4$ (right). Note that $\Gamma=\Gamma^x\cup \Gamma^y$, all $\Gamma^x_p$ and $\Gamma^y_q$ are open sets for $1\le p,q \le 3$. The intersection points $\overline{\Gamma^x}\cap \overline{\Gamma^y}$ are indicated by the blue color and points on $\overline{\Gamma}\cap \partial \Omega$ are indicated by the black color.}
		\label{Gamma}
\end{figure}

In this paper,  we discuss the following elliptic cross-interface problem (see \cref{coefficient,Gamma} for illustrations): {\em Let the domain $\Omega:=(0,1)^2$. Then we consider:
\be \label{intersect:3}
	\left\{ \begin{array}{llll}
		-\nabla \cdot \Big(a \nabla  u\Big)&=&f 
&\mbox{in } \Omega \setminus {\Gamma},\\
		\left[u\right]&=&0 
&\mbox{on }\Gamma,\\	
		\left[a \nabla  u \cdot \vec{n}\right] &=&0 
&\mbox{on }\Gamma,\\
		u&=&0 
&\mbox{on }	\partial\Omega,
	\end{array}
	\right.
\ee
where $\vec{n}$ is the unit normal vector of  $\Gamma$ and
\be\label{GammaXY}
\begin{split}
& \Gamma:=\Gamma^x \cup \Gamma^y, \qquad \Gamma^x:=\displaystyle\cup_{p=1}^{m-1}  \Gamma^x_p,\qquad \Gamma^y:=\displaystyle\cup_{q=1}^{m-1}   \Gamma^y_q, \\
& \Gamma^x_p:= \Big\{(x,y) : x=\frac{p}{m},\ y\in (0,1), \text{ and } \bmod(y m,1)\ne 0\Big \}, \qquad 1\le p \le m-1,\\
& \Gamma^y_q:= \Big\{(x,y) : y=\frac{q}{m},\ x\in (0,1), \text{ and } \bmod(x m,1)\ne 0\Big \}, \qquad 1\le q \le m-1,
\end{split}
\ee
where $m$ is a positive integer.} Note that 
\[
\Omega,\ \Gamma^x_p,\ \Gamma^y_q,\ \Gamma^x,\ \Gamma^y,\ \Gamma,\ \Omega \setminus \overline{\Gamma} \text{ are all open sets for every } 1\le p,q \le m-1.
\]
For $(\xi,y)\in \Gamma^x$ (i.e., on the vertical line of the cross-interface $\Gamma$),
\[
[u](\xi,y):=\lim_{x \to \xi^+}u(x,y)- \lim_{ x\to \xi^-}u(x,y),
\quad
[a \nabla  u \cdot \vec{n}](\xi,y):=\lim_{x \to \xi^+} a(x,y) \frac{\partial u}{\partial x}(x,y) - \lim_{x\to \xi^-} a(x,y)\frac{\partial u}{\partial x}(x,y);
\]
while for $(x,\zeta)\in \Gamma^y$  (i.e., on the horizontal line of the cross-interface $\Gamma$),
\[
[u](x,\zeta):=\lim_{y\to \zeta^+ }u(x,y)- \lim_{ y\to \zeta^-}u(x,y),  \quad
[a \nabla  u \cdot \vec{n}](x,\zeta):=\lim_{y \to \zeta^+} a(x,y) \frac{\partial u}{\partial y}(x,y) - \lim_{y\to \zeta^-} a(x,y)\frac{\partial u}{\partial y}(x,y).
\]
We also define the subdomains of $\Omega$ as follows:
\be\label{Def:Omegaij}
\Omega_{p,q}:=\Big\{(x,y) : \frac{p-1}{m}<x<\frac{p}{m},\ \frac{q-1}{m}<y<\frac{q}{m} \Big\}, \qquad 1\le p,q \le m.
\ee
Clearly, $a$ is a constant in each $\Omega_{p,q}$ and
\[
\cup_{p=1}^m\cup_{q=1}^m \Omega_{p,q}=\Omega\setminus \overline{\Gamma}.
\]

In this paper, we discuss the model problem \eqref{intersect:3} under the following assumptions:
\begin{itemize}
	
	\item[(A1)]  The coefficient function $a$ is a positive piecewise  constant function in $\Omega \setminus {\Gamma}$ (see \cref{coefficient} for an illustration).
	
	\item[(A2)]  The source term $f$ is continuous in $\Omega$.
\end{itemize}
\begin{figure}[htbp]
	\centering	
	\includegraphics[width=9cm,height=9cm]{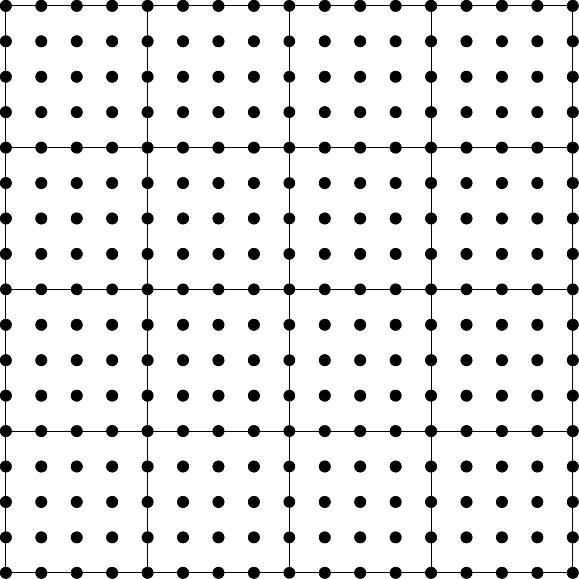}
	\caption{An example for uniform Cartesian grids of the problem \eqref{intersect:3}}
	\label{uniform:mesh}
\end{figure}

The rest of this paper is structured as follows:

In \cref{section:FEMandFDM}, we  present  second-order finite element and finite difference  methods using uniform Cartesian grids. Precisely, we illustrate  second-order finite element method in \cref{Sec:FEM} and second-order finite difference method in \cref{Sec:FDM}.

In \cref{Intersect:sec:Numeri}, we provide 7 numerical examples in the following two cases:
\begin{itemize}
	\item We examine the model problem \eqref{intersect:3} using a high-contrast and low-frequency coefficient function in \cref{Intersect:ex1}, and a low-contrast and relatively high-frequency coefficient function in \cref{Intersect:ex2,Intersect:ex3:New}. Based on performances of the FEM and FDM solutions presented in \cref{fig:exam1,fig:exam2,fig:exam3:New}, we observe that the FEM and FDM yield similar numerical solutions.
	\item We test the model problem \eqref{intersect:3} with high-contrast and relatively high-frequency coefficient functions  in \cref{Intersect:ex3,Intersect:ex4,Intersect:ex5,Intersect:ex7}. From  performances of the FEM and FDM solutions in \cref{fig:exam3,fig:exam4,fig:exam5,fig:exam7}, we observe that FEM and FDM produce completely different numerical solutions.
\end{itemize}

In \cref{Intersect:sec:Reason}, we present several possible reasons to explain the  above numerical results.

In \cref{sec:Conclu}, we highlight the key contributions of this paper.

\section{Second-order finite element and finite difference  methods using uniform Cartesian grids}
\label{section:FEMandFDM}
In this section, we present second-order finite element and finite difference  methods on uniform Cartesian grids for the elliptic cross-interface problem in \eqref{intersect:3}. 
We  define that
\be\label{xiyj}
x_i:=i h, \quad i=0,1,\ldots,N, \quad \text{and} \quad y_j:=j h, \quad j=0,1,\ldots,N,
\ee
where the uniform Cartesian mesh size $h:=\frac{1}{N}$ is chosen for some positive integer $N$ such that all intersection points $(x,y)\in \overline{\Gamma^x}\cap \overline{\Gamma^y}$ are grid points and some grid points lie on the vertical/horizontal interface lines (see \cref{uniform:mesh} for an example with $m=4$, $h=1/16$ and $N=16$).
\subsection{Second-order finite element method}
\label{Sec:FEM}
We use the standard second-order finite element method with  piecewise linear basis and continuous functions to solve \eqref{intersect:3}. Precisely,
let's define the following 1D piecewise linear and continuous functions on $\{x_0,x_1,\dots,x_{N}\}$ and  $\{y_0,y_1,\dots,y_{N}\}$ respectively:
\[
\phi_i(x)
:=\begin{cases}
	\displaystyle\frac{x-x_{i-1}}{x_i-x_{i-1}} & \text{if } x_{i-1} \le x \le x_{i} ,\\
	\displaystyle\frac{x_{i+1}-x}{x_{i+1}-x_i} & \text{if } x_{i} \le x \le x_{i+1} ,\\
	0  & \text{else},
\end{cases}
\qquad
\phi_j(y)
:=\begin{cases}
	\displaystyle\frac{y-y_{j-1}}{y_j-y_{j-1}} & \text{if } y_{j-1} \le y \le y_{j} ,\\
	\displaystyle\frac{y_{j+1}-y}{y_{j+1}-y_j} & \text{if } y_{j} \le y \le y_{j+1} ,\\
	0  & \text{else},
\end{cases}
\]
where  $1\le i,j\le N-1$.
Then we define the  2D basis functions, and the continuous piecewise linear finite element space as follows:
\[
\phi_{i,j}:=\phi_i(x)\phi_j(y), \quad 1\le i,j\le N-1,
\]
\[
V_h:=\{ \phi_{i,j} : 1\le i,j\le N-1 \}.
\]
Let $u^E_h$ be the  numerical finite element method solution of the exact solution $u$ of the elliptic cross-interface problem \eqref{intersect:3}  using the mesh size $h=1/N$.
Then the standard finite element method implies
\[
u^E_h:=\sum_{i=1}^{N-1}\sum_{j=1}^{N-1} c_{i,j}\phi_{i,j},
\]
where every $c_{i,j}$ is to-be-determined constant for $1\le i,j\le N-1$ by solving the following weak formulation of \eqref{intersect:3}:
\[
\int_\Omega a \nabla u^E_h\cdot \nabla v=\int_\Omega f v,\ \ \ \ \mbox{for all}\ v\in V_h.
\]
To verify the convergence rate of the finite element solution, we also define that
\be 
(u^E_h)_{i,j}=u^E_h(x_i,y_j), \quad \text{for} \quad 1\le i,j\le N-1,
\ee
where $(x_i,y_j)$ is defined in $\eqref{xiyj}$.
\subsection{Second-order finite difference method}
\label{Sec:FDM}
\begin{figure}[tbhp]
	\centering	
	\begin{subfigure}[b]{0.3\textwidth}
	\includegraphics[width=5cm,height=5cm]{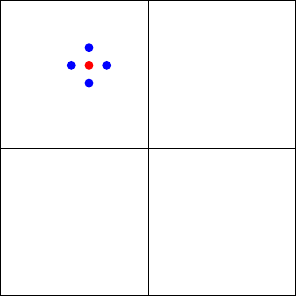}
\end{subfigure}
\begin{subfigure}[b]{0.3\textwidth}
	\includegraphics[width=5cm,height=5cm]{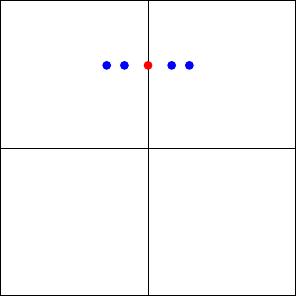}
\end{subfigure}
\begin{subfigure}[b]{0.3\textwidth}
	\includegraphics[width=5cm,height=5cm]{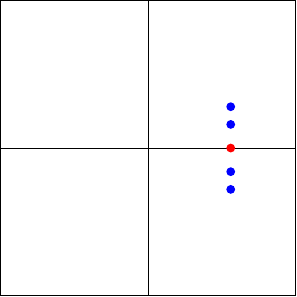}
\end{subfigure}
	\caption{ First panel:  5-point  finite difference scheme in \cref{thm:regular} with $(x_i,y_j)\in \Omega\setminus \Gamma$. Second panel:   5-point  finite difference scheme in \cref{thm:irregular:V} with $(x_i,y_j)\in \Gamma^x$. Third panel:  5-point  finite difference scheme in \cref{thm:irregular:H} with $(x_i,y_j)\in \Gamma^y$. The grid point $(x_i,y_j)$ is marked in red.
	}
	\label{fig:interface}
\end{figure}

We use the standard second-order 5-point finite difference method for the interior grid point $(x_i,y_j)\in \Omega \setminus \Gamma$ (see the first panel of \cref{fig:interface}). For the interface grid point $(x_i,y_j)\in  \Gamma:=\Gamma^x \cup \Gamma^y$ (see the second and third panels of \cref{fig:interface}), we use second-order 5-point finite difference methods based on second-order 3-point backward/forward upwind schemes.

Let $(u^D_h)_{i,j}$ be the value of the numerical finite difference method solution $u^D_h$ of the exact solution $u$ of the elliptic cross-interface problem \eqref{intersect:3}, at the grid point $(x_i, y_j)$ using the mesh size $h=1/N$.  For the sake of brevity,
 we also define that
\[
u_{i,j}=u(x_i,y_j),\qquad a_{i,j}=a(x_i,y_j),\qquad  f_{i,j}=f(x_i,y_j).
\]
Then we use the well-known 5-point finite difference method at $(x_i,y_j)\in \Omega \setminus \Gamma$ in the following theorem:
\begin{theorem}\label{thm:regular}
	Consider the interior grid point $(x_i,y_j)\in \Omega \setminus \Gamma$. Then the following 5-point finite difference scheme  (see the first panel of \cref{fig:interface})
\[
a_{i,j}\bigg(\big(u^D_h\big)_{i-1,j}+\big(u^D_h\big)_{i+1,j}+\big(u^D_h\big)_{i,j-1}+\big(u^D_h\big)_{i,j+1}-4\big(u^D_h\big)_{i,j}\bigg)=-h^2f_{i,j},
\]
is the second-order consistent at the interior grid point $(x_i,y_j)\in \Omega \setminus \Gamma$.
\end{theorem}
For the interface grid point $(x_i,y_j)\in  \Gamma^x$,
we have that the
second-order 3-point backward upwind scheme to approximate $u_x(x_i,y_j)$
is
\be 
\frac{\partial u}{\partial x}(x_i,y_j)=\frac{3u_{i,j}-4u_{i-1,j}+u_{i-2,j}}{2h}+\bo(h^{2}),
\ee
and
the second-order 3-point forward upwind scheme to approximate $u_x(x_i,y_j)$ is
\be 
\frac{\partial u}{\partial x}(x_i,y_j)=\frac{-3u_{i,j}+4u_{i+1,j}-u_{i+2,j}}{2h}+\bo(h^{2}).
\ee
Since  the coefficient function $a$ is a piecewise  constant function in $\Omega \setminus {\Gamma}$ (see \cref{coefficient,uniform:mesh} for illustrations), we have 
\be
\begin{split}
[a \nabla  u \cdot \vec{n}](x_i,y_j)=&a(x_i+h,y_j)\frac{-3u_{i,j}+4u_{i+1,j}-u_{i+2,j}}{2h}\\
&-a(x_i-h,y_j)\frac{3u_{i,j}-4u_{i-1,j}+u_{i-2,j}}{2h}+\bo(h^{2}),
\end{split}
\ee
for $(x_i,y_j)\in \Gamma^x$. By $[a \nabla  u \cdot \vec{n}]=0$ on $\Gamma^x$, we obtain that
\[
\frac{1}{2h}\Big(a(x_i+h,y_j)\big(-3u_{i,j}+4u_{i+1,j}-u_{i+2,j}\big)-a(x_i-h,y_j)\big(3u_{i,j}-4u_{i-1,j}+u_{i-2,j}\big)\Big)=0,
\]
	has a second order of consistency at the interface grid point $(x_i,y_j)\in  \Gamma^x$.
	
\begin{figure}[htbp]
	\centering	
	\includegraphics[width=9cm,height=9cm]{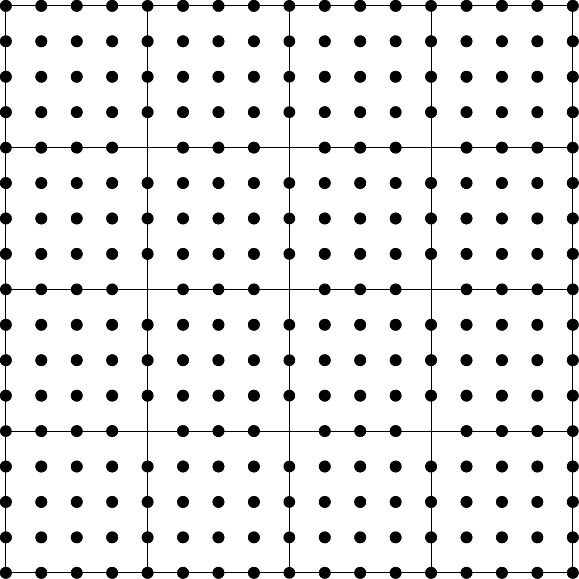}
	\caption{The computed numerical solutions $\big(u^D_h\big)_{i,j}$ with  $0\le i,j \le N$, and  $(x_i,y_j)\notin  \overline{\Gamma^x}\cap \overline{\Gamma^y}$
in the finite difference method. 	
}
	\label{ignore:intersection}
\end{figure}

	In summary, we have the second-order 5-point finite difference schemes at $(x_i,y_j)\in  \Gamma:=\Gamma^x\cup \Gamma^y$ in the following \cref{thm:irregular:V,thm:irregular:H}.
\begin{theorem}\label{thm:irregular:V}
	Consider the interface grid point $(x_i,y_j)\in \Gamma^x$. Let $a^{-}=a(x_i-h,y_j)$ and $a^{+}=a(x_i+h,y_j)$. Then the following 5-point finite difference scheme  (see the second panel of \cref{fig:interface})
	\[
-a^-\big(u^D_h\big)_{i-2,j}+4a^-\big(u^D_h\big)_{i-1,j}-3(a^-+a^+)\big(u^D_h\big)_{i,j}+4a^+\big(u^D_h\big)_{i+1,j}-a^+\big(u^D_h\big)_{i+2,j}=0,
	\]
	has a second order of consistency at the interface grid point $(x_i,y_j)\in  \Gamma^x$.
\end{theorem}
\begin{theorem}\label{thm:irregular:H}
	Consider the interface grid point $(x_i,y_j)\in \Gamma^y$. Let $a^{-}=a(x_i,y_j-h)$ and $a^{+}=a(x_i,y_j+h)$. Then  the following 5-point finite difference scheme  (see the third panel of \cref{fig:interface})
	\[
-a^-\big(u^D_h\big)_{i,j-2}+4a^-\big(u^D_h\big)_{i,j-1}-3(a^-+a^+)\big(u^D_h\big)_{i,j}+4a^+\big(u^D_h\big)_{i,j+1}-a^+\big(u^D_h\big)_{i,j+2}=0,
\]
	has a second order of consistency at the interface grid point $(x_i,y_j)\in  \Gamma^y$.
\end{theorem}
\begin{remark}
The second-order 5-point schemes in  \cref{thm:irregular:V,thm:irregular:H} can also be derived from \cite[(A.6)-(A.9)]{FengHanMinevCross} by choosing $M=2$ and applying two homogeneous jump conditions.
\end{remark}
By the definition of the uniform Cartesian grid in \eqref{xiyj} of the model problem \eqref{intersect:3}, we also need to provide the finite difference scheme at the intersection grid point $(x_i,y_j)\in  \overline{\Gamma^x}\cap \overline{\Gamma^y}$ (see \cref{uniform:mesh} for an illustration). But the schemes in \cref{thm:regular,thm:irregular:V,thm:irregular:H} do not use any intersection grid points $(x_i,y_j)\in  \overline{\Gamma^x}\cap \overline{\Gamma^y}$. So it is not necessary to propose the finite difference scheme at the grid point  $\{ (x_i,y_j) : (x_i,y_j)\in  \overline{\Gamma^x}\cap \overline{\Gamma^y}\}$. I.e., we only compute 
\[
\big(u^D_h\big)_{i,j}, \quad 0\le i,j \le N,\quad  \text{and} \quad (x_i,y_j)\notin  \overline{\Gamma^x}\cap \overline{\Gamma^y},
\]
in the finite difference method (see \cref{ignore:intersection} for an illustration). 
\section{Numerical experiments}\label{Intersect:sec:Numeri}
Recall that $\big(u^E_{h}\big)_{i,j}$ and $\big(u^D_{h}\big)_{i,j}$ are the finite element and finite difference solutions $u^E_{h}$ and $u^D_{h}$ respectively  of the model problem \eqref{intersect:3} at $(x_i, y_j)$ using the mesh size $h$, where
\[
x_i:=i h, \quad i=0,1,\ldots,N, \quad \text{and} \quad y_j:=j h, \quad j=0,1,\ldots,N,
\]
$h:=1/N$ for some integer $N$, and $\Omega:=(0,1)^2$.

\begin{figure}[htbp]
	\centering
	\begin{subfigure}[b]{0.4\textwidth}
		\hspace{-1cm}
	\includegraphics[width=7cm,height=7cm]{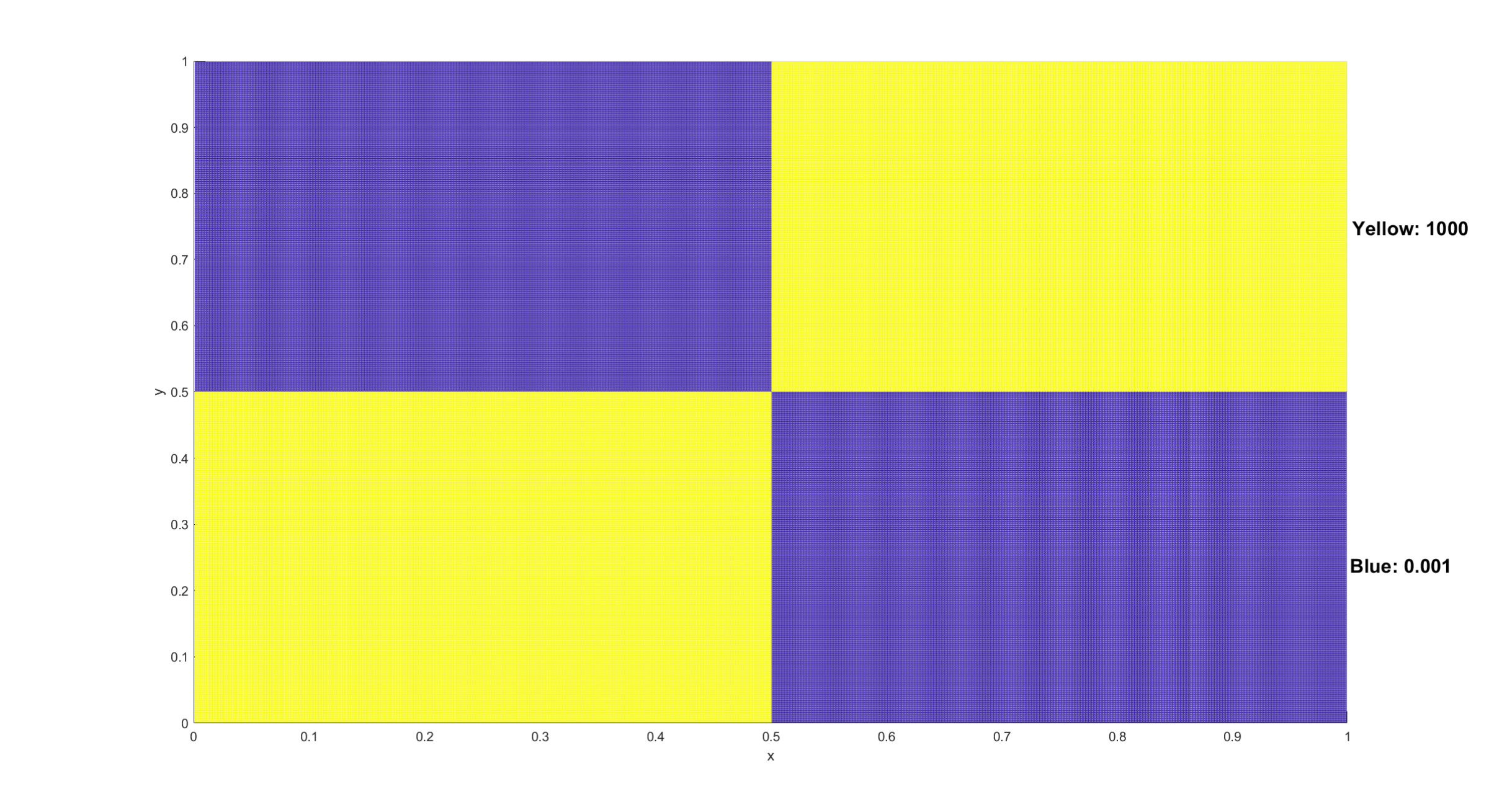}
\end{subfigure}
	\begin{subfigure}[b]{0.4\textwidth}
		\hspace{0.7cm}	
		\vspace{0.8cm}	
	\includegraphics[width=5.6cm,height=5.6cm]{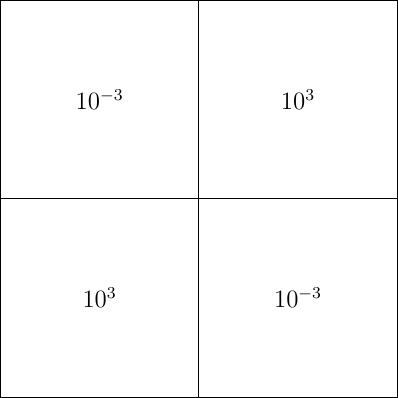}
\end{subfigure}
	\begin{subfigure}[b]{0.45\textwidth}
		\includegraphics[width=7cm,height=7cm]{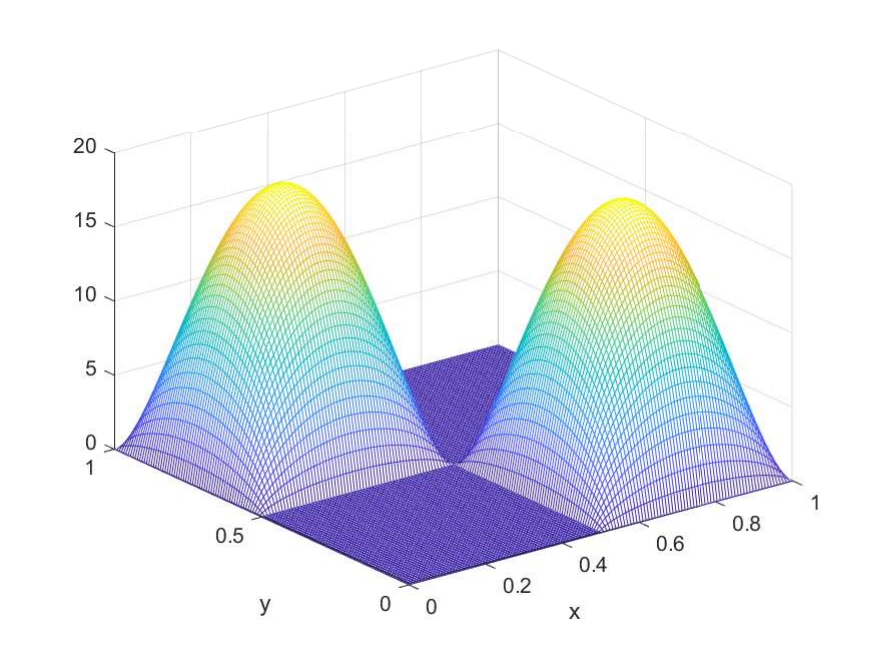}
	\end{subfigure}
	\begin{subfigure}[b]{0.45\textwidth}
		\includegraphics[width=7cm,height=7cm]{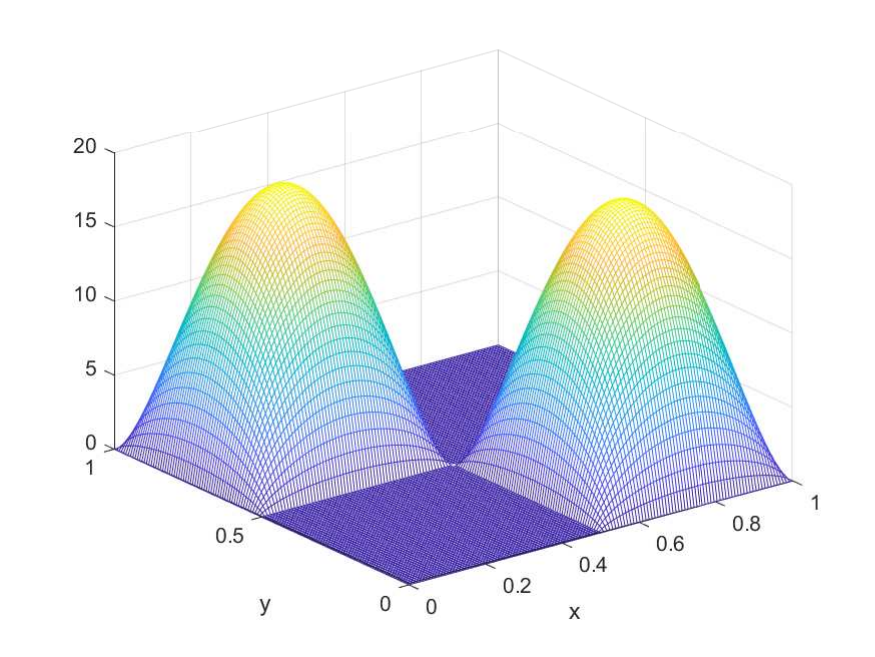}
	\end{subfigure}
	\begin{subfigure}[b]{0.45\textwidth}
		\includegraphics[width=7cm,height=7cm]{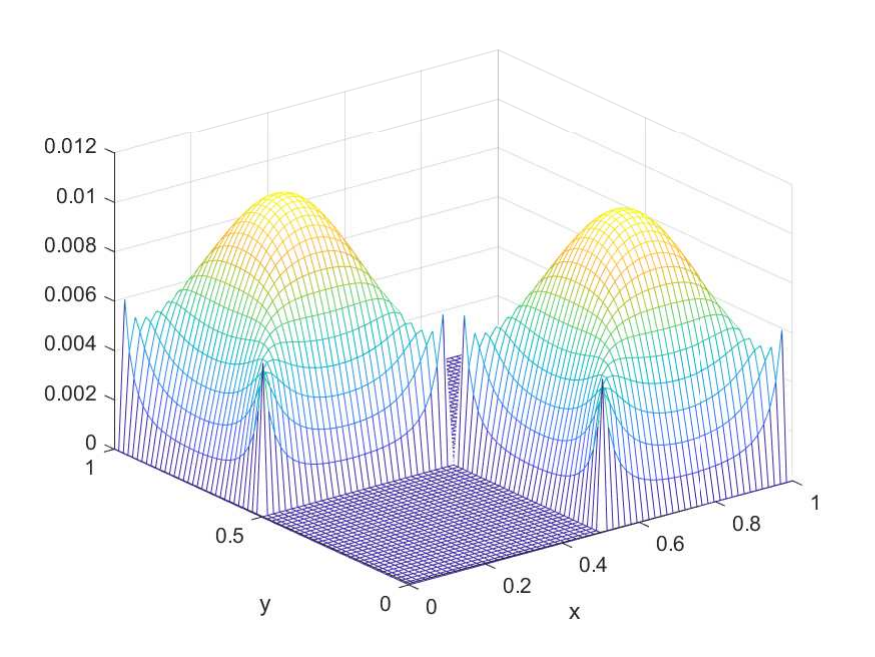}
	\end{subfigure}
	\begin{subfigure}[b]{0.45\textwidth}
		\includegraphics[width=7cm,height=7cm]{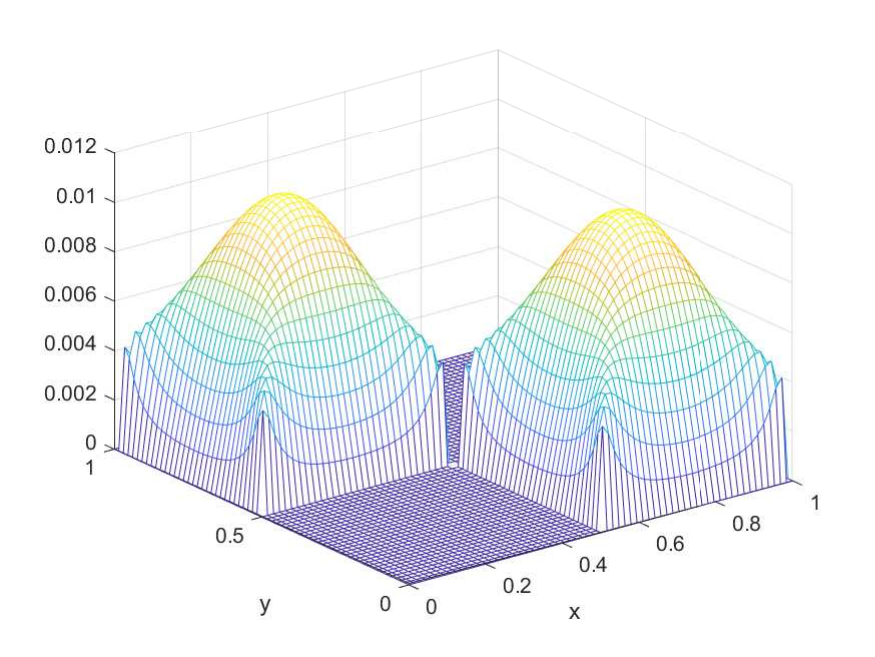}
	\end{subfigure}
	\caption{ \cref{Intersect:ex1}: The first row: the coefficient function $a$ with $a=1000$ in yellow squares, $a=0.001$ in blue squares, and $m=2$. The second row: the finite element solution $(u^E_h)_{i,j}$ (left) and the finite difference solution $(u^D_h)_{i,j}$ (right) at all  $(x_i,y_j)$ with $h=1/2^7$. The third row: $\big|(u^E_{h})_{i,j}-(u^E_{h/2})_{2i,2j}\big|$ (left) and  $\big|(u^D_{h})_{i,j}-(u^D_{h/2})_{2i,2j}\big|$ (right) at all $(x_i,y_j)$ with $h=1/2^6$.}
	\label{fig:exam1}
\end{figure}

\begin{figure}[htbp]
	\centering
	\begin{subfigure}[b]{0.4\textwidth}
	\hspace{-1cm}
	\includegraphics[width=7cm,height=7cm]{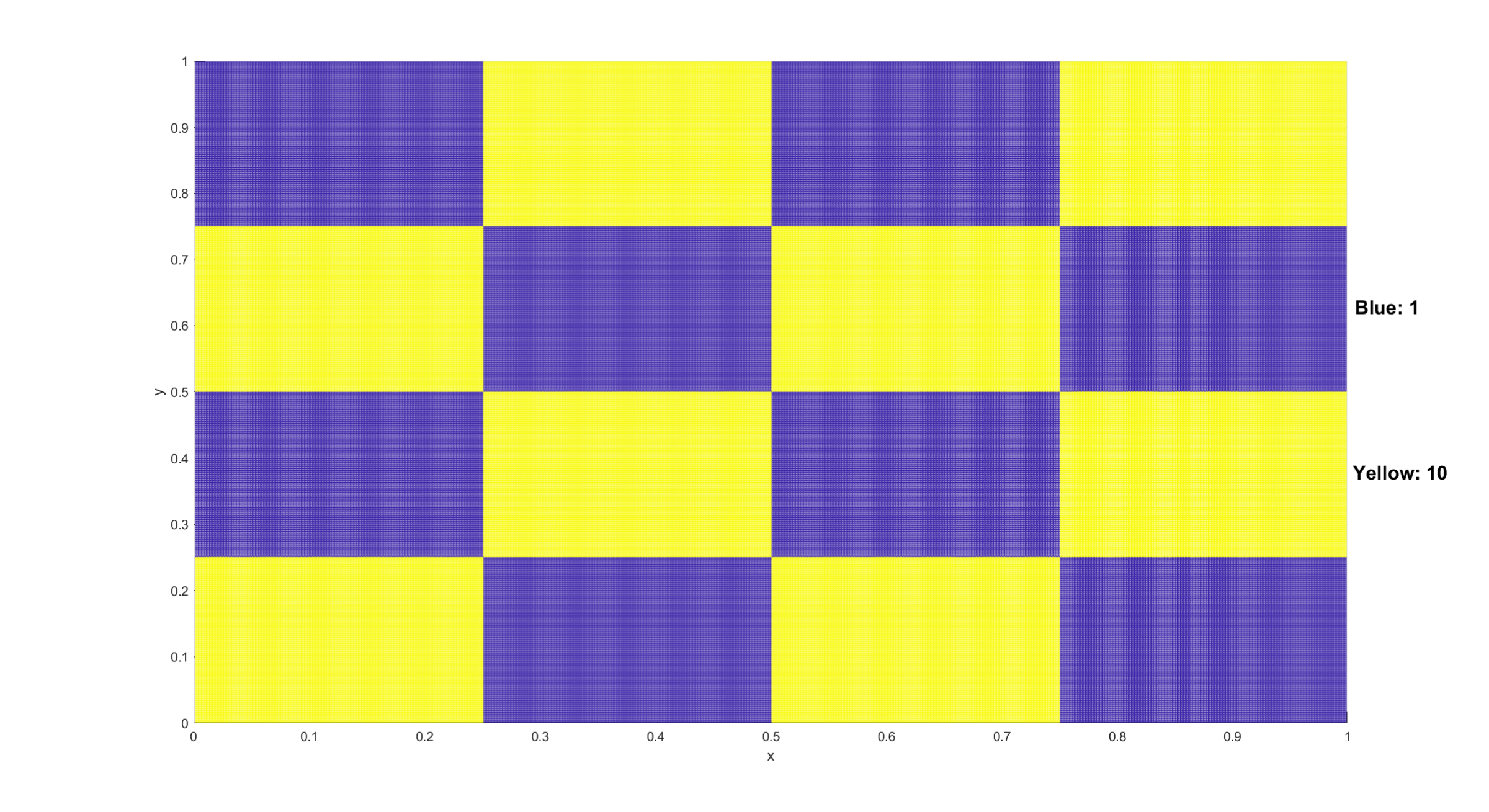}
\end{subfigure}
\begin{subfigure}[b]{0.4\textwidth}
	\hspace{0.7cm}	
	\vspace{0.8cm}	
	\includegraphics[width=5.6cm,height=5.6cm]{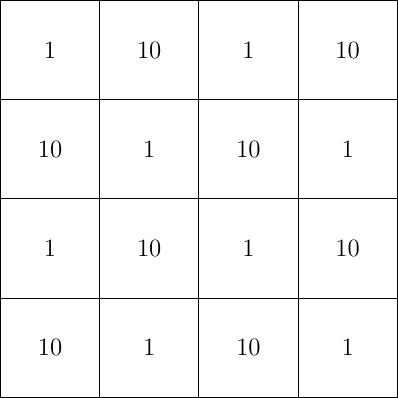}
\end{subfigure}
	\begin{subfigure}[b]{0.45\textwidth}
		\includegraphics[width=7cm,height=7cm]{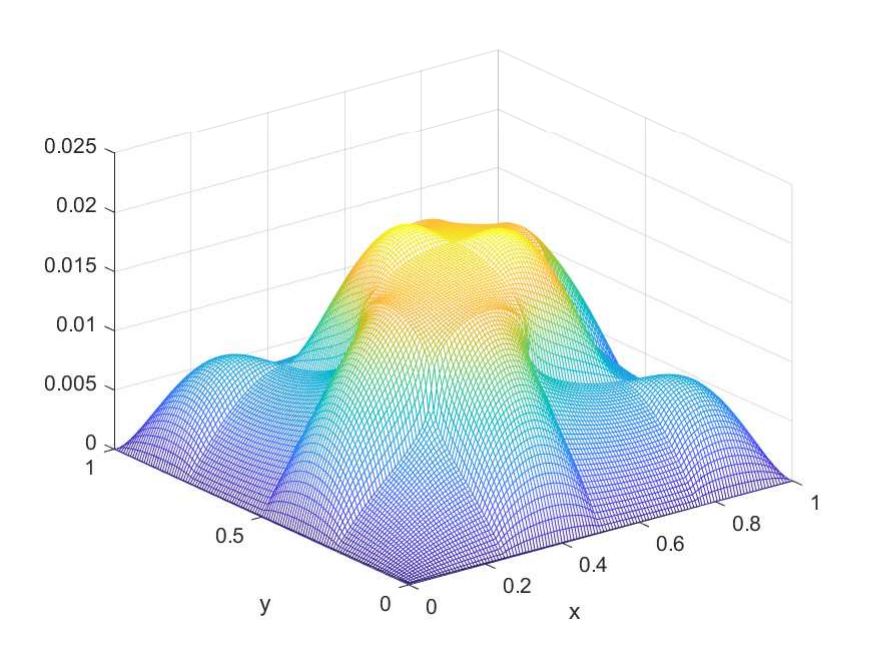}
	\end{subfigure}
	\begin{subfigure}[b]{0.45\textwidth}
		\includegraphics[width=7cm,height=7cm]{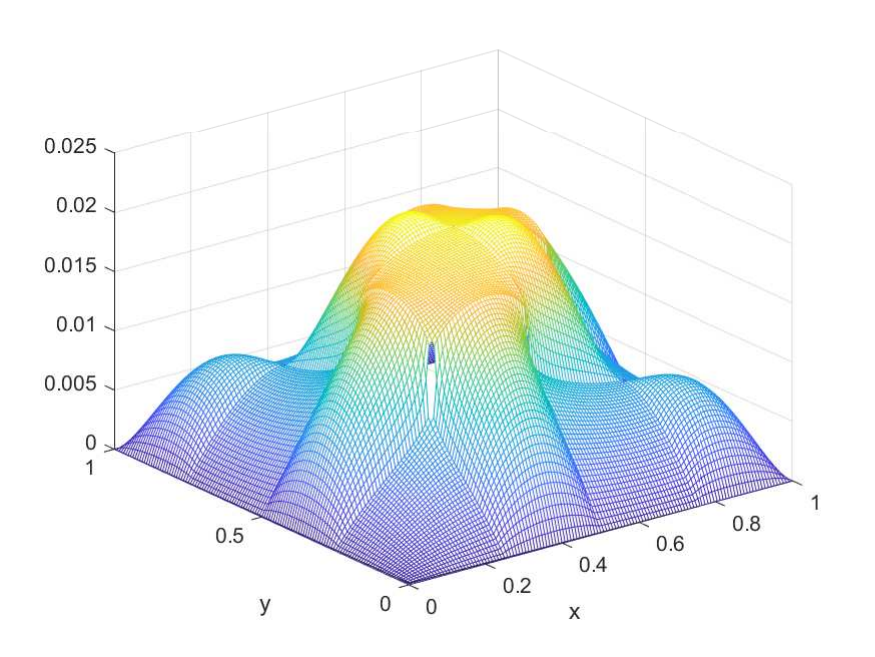}
	\end{subfigure}
	\begin{subfigure}[b]{0.45\textwidth}
		\includegraphics[width=7cm,height=7cm]{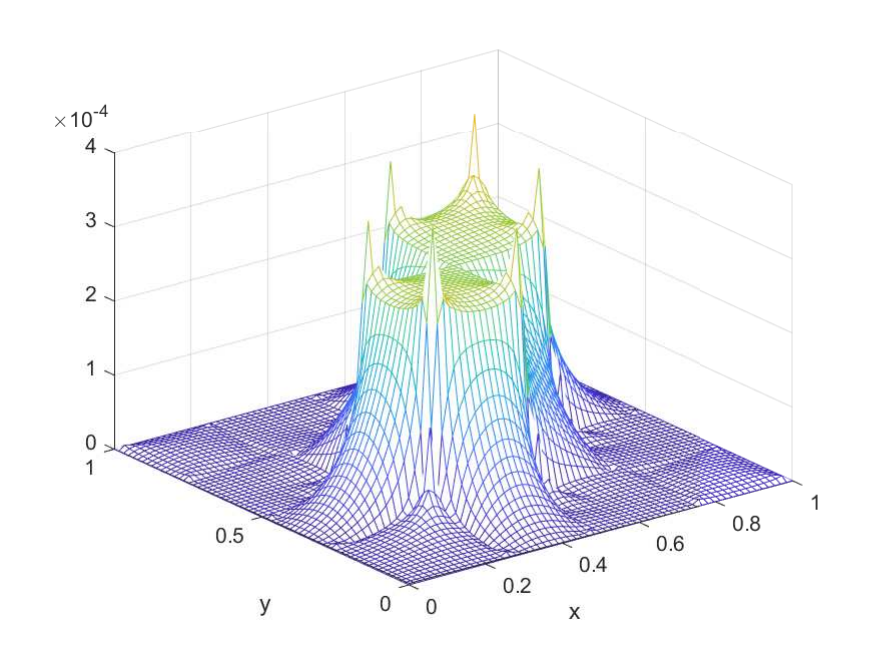}
	\end{subfigure}
	\begin{subfigure}[b]{0.45\textwidth}
		\includegraphics[width=7cm,height=7cm]{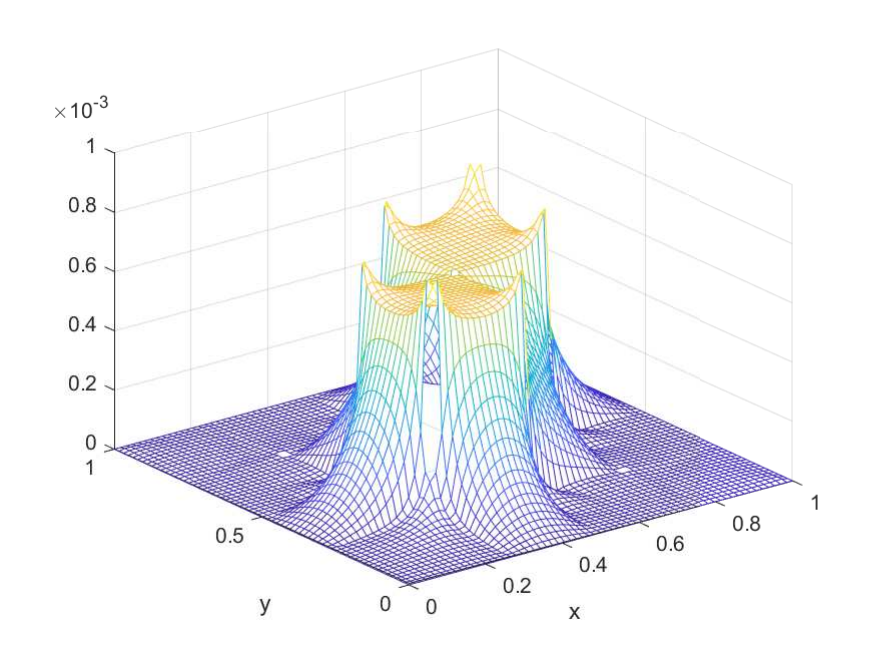}
	\end{subfigure}
	\caption{ \cref{Intersect:ex2}: The first row: the coefficient function $a$ with $a=10$ in yellow squares, $a=1$ in blue squares, and $m=4$. The second row: the finite element solution $(u^E_h)_{i,j}$ (left) and the finite difference solution $(u^D_h)_{i,j}$ (right) at all $(x_i,y_j)$ with $h=1/2^7$. The third row: $\big|(u^E_{h})_{i,j}-(u^E_{h/2})_{2i,2j}\big|$ (left) and  $\big|(u^D_{h})_{i,j}-(u^D_{h/2})_{2i,2j}\big|$ (right) at all  $(x_i,y_j)$ with $h=1/2^6$.}
\label{fig:exam2}
\end{figure}

\begin{figure}[htbp]
	\centering
	\begin{subfigure}[b]{0.4\textwidth}
	\hspace{-1cm}
	\includegraphics[width=7cm,height=7cm]{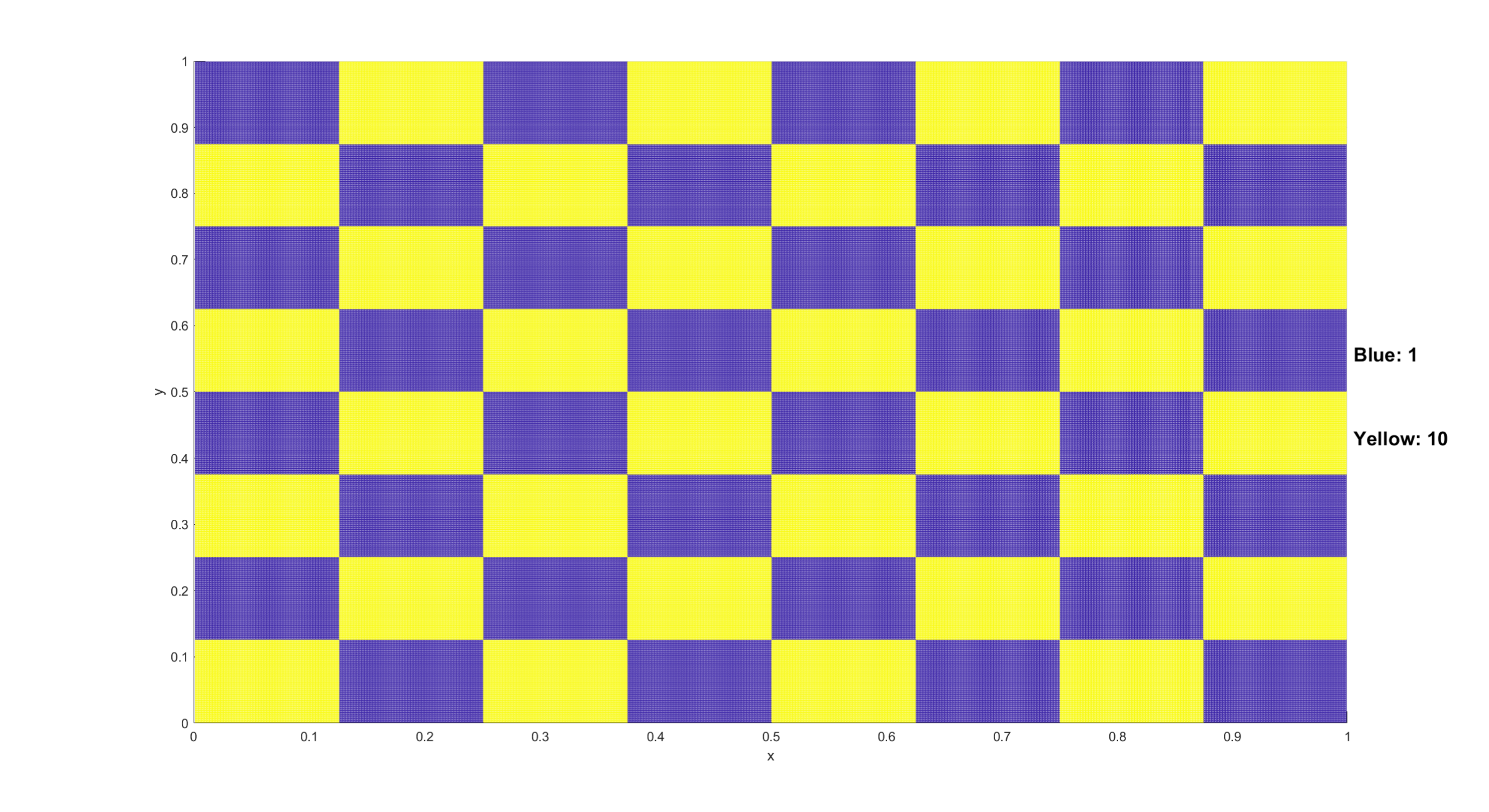}
\end{subfigure}
\begin{subfigure}[b]{0.4\textwidth}
	\hspace{0.7cm}	
	\vspace{0.8cm}	
	\includegraphics[width=5.6cm,height=5.6cm]{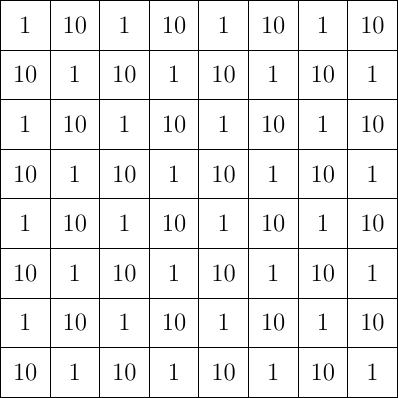}
\end{subfigure}
	\begin{subfigure}[b]{0.45\textwidth}
		\includegraphics[width=7cm,height=7cm]{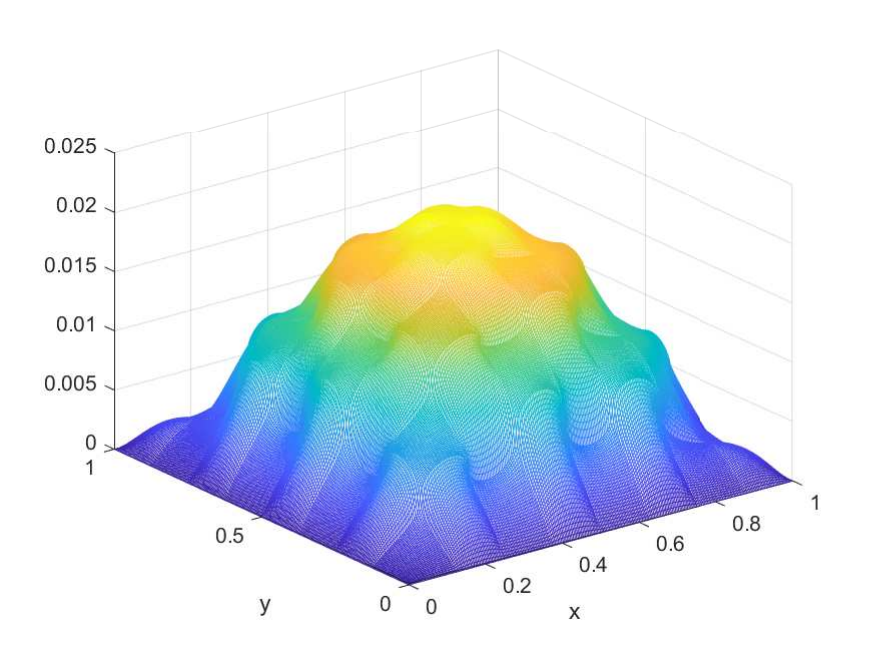}
	\end{subfigure}
	\begin{subfigure}[b]{0.45\textwidth}
		\includegraphics[width=7cm,height=7cm]{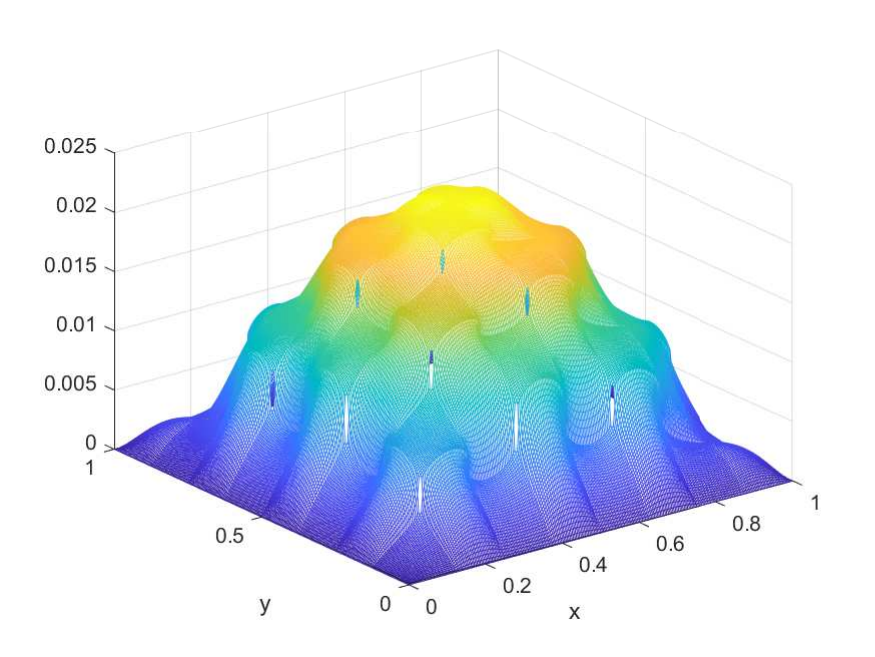}
	\end{subfigure}
	\begin{subfigure}[b]{0.45\textwidth}
		\includegraphics[width=7cm,height=7cm]{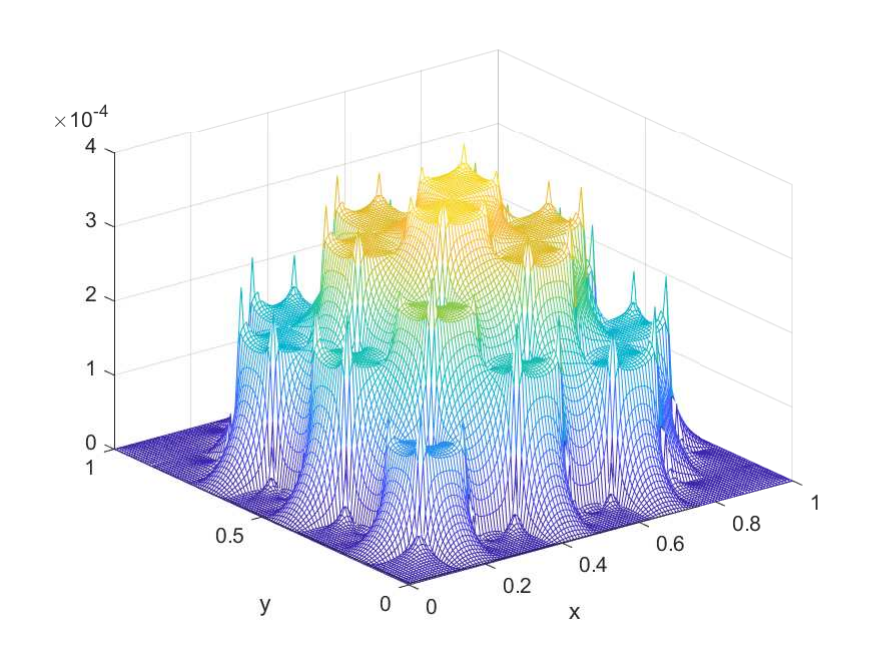}
	\end{subfigure}
	\begin{subfigure}[b]{0.45\textwidth}
		\includegraphics[width=7cm,height=7cm]{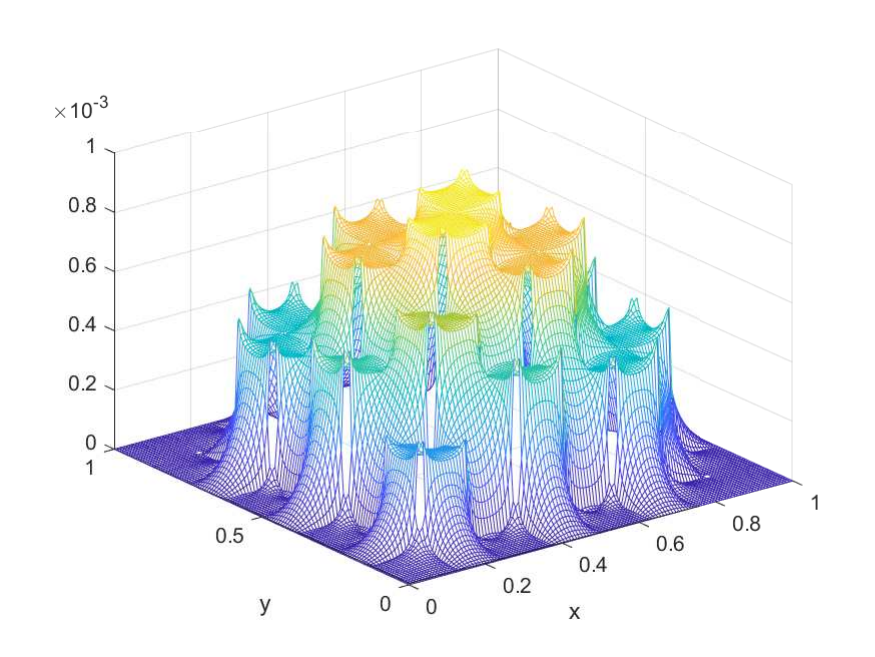}
	\end{subfigure}
	\caption{ \cref{Intersect:ex3:New}: The first row: the coefficient function $a$ with $a=10$ in yellow squares, $a=1$ in blue squares, and $m=8$. The second row: the finite element solution $(u^E_h)_{i,j}$ (left) and the finite difference solution $(u^D_h)_{i,j}$ (right) at all  $(x_i,y_j)$ with $h=1/2^8$. The third row: $\big|(u^E_{h})_{i,j}-(u^E_{h/2})_{2i,2j}\big|$ (left) and  $\big|(u^D_{h})_{i,j}-(u^D_{h/2})_{2i,2j}\big|$ (right) at all  $(x_i,y_j)$ with $h=1/2^7$.}
	\label{fig:exam3:New}
\end{figure}

\begin{figure}[htbp]
	\centering
	\begin{subfigure}[b]{0.4\textwidth}
	\hspace{-1cm}
	\includegraphics[width=7cm,height=7cm]{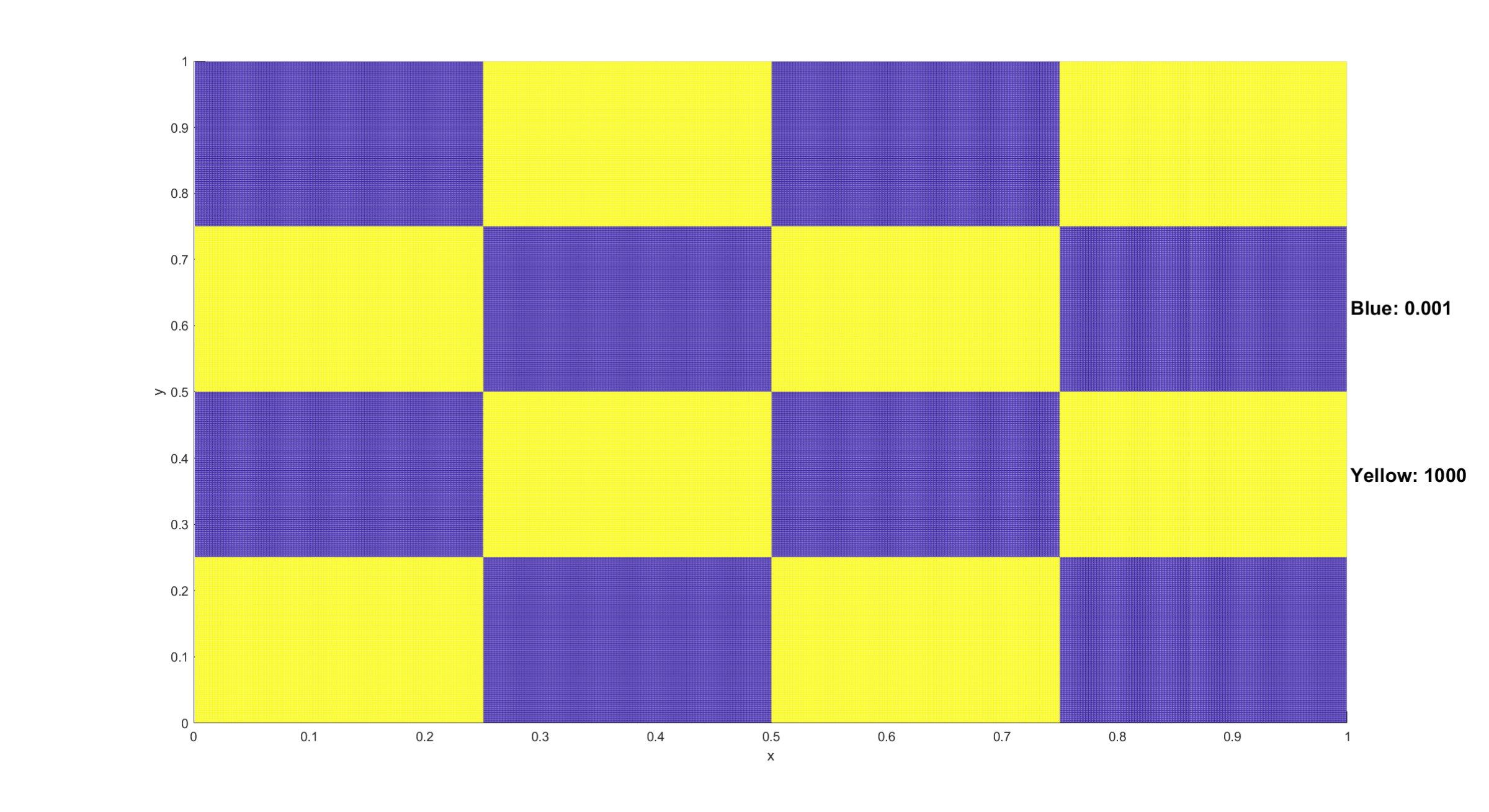}
\end{subfigure}
\begin{subfigure}[b]{0.4\textwidth}
	\hspace{0.7cm}	
	\vspace{0.8cm}	
	\includegraphics[width=5.6cm,height=5.6cm]{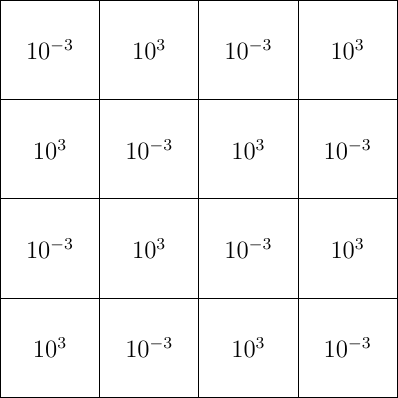}
\end{subfigure}
	\begin{subfigure}[b]{0.45\textwidth}
		\includegraphics[width=7cm,height=7cm]{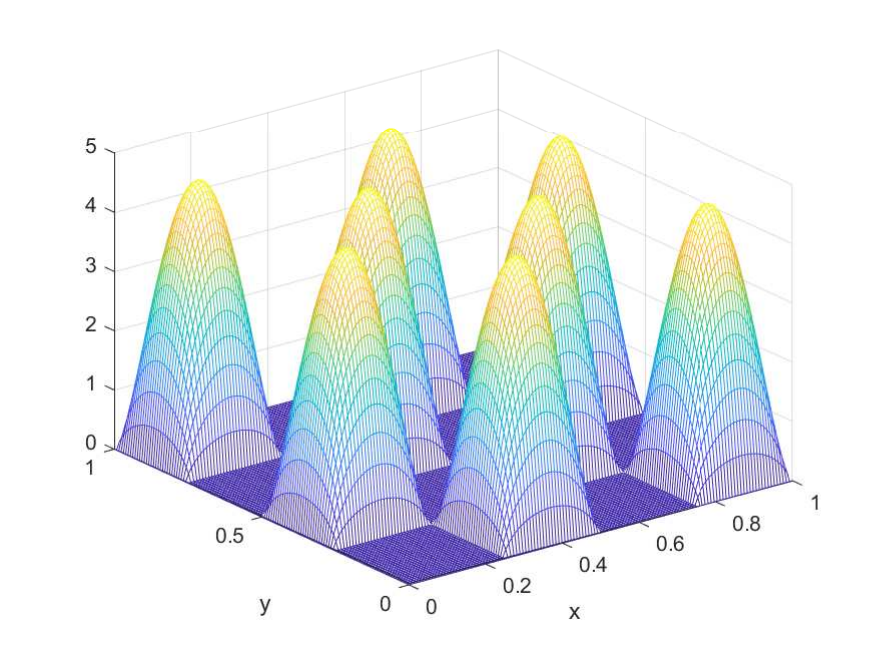}
	\end{subfigure}
	\begin{subfigure}[b]{0.45\textwidth}
		\includegraphics[width=7cm,height=7cm]{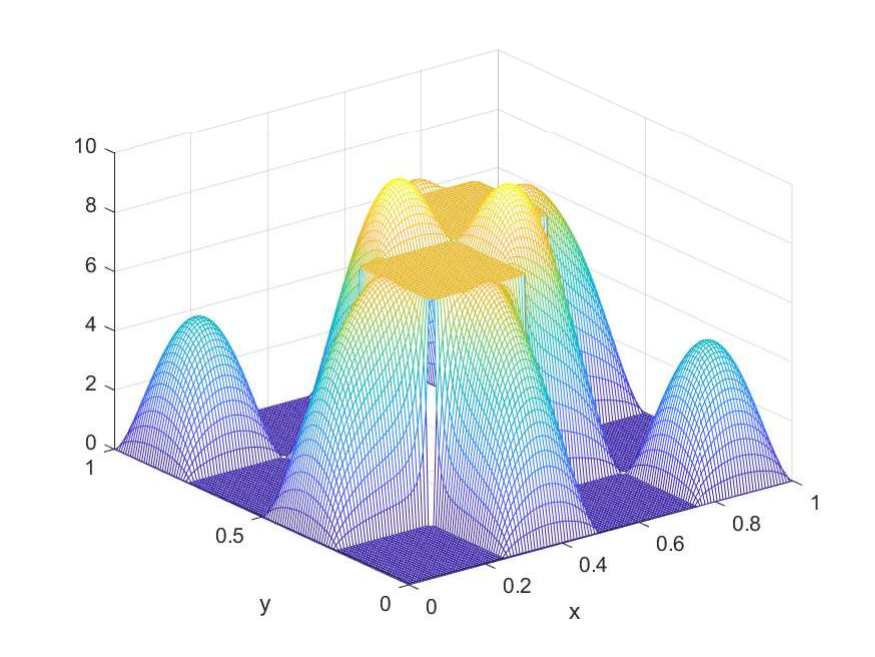}
	\end{subfigure}
	\begin{subfigure}[b]{0.45\textwidth}
		\includegraphics[width=7cm,height=7cm]{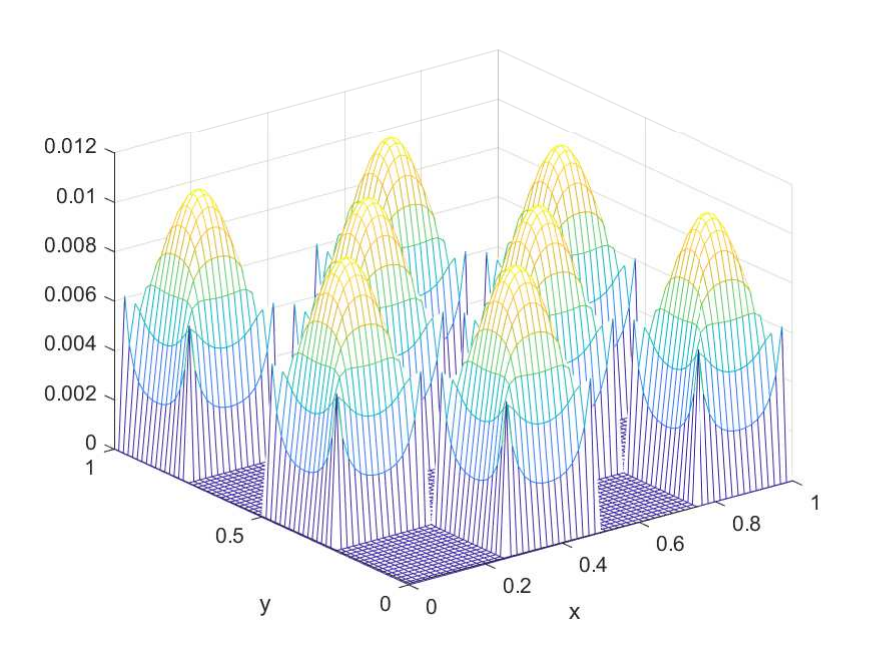}
	\end{subfigure}
	\begin{subfigure}[b]{0.45\textwidth}
		\includegraphics[width=7cm,height=7cm]{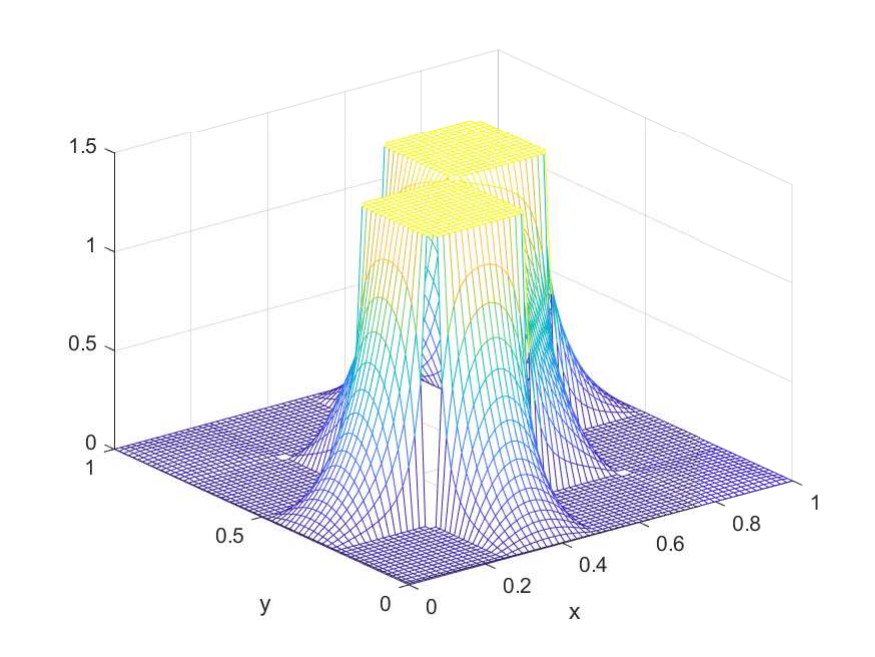}
	\end{subfigure}
	\caption{ \cref{Intersect:ex3}: The first row: the coefficient function $a$ with $a=1000$ in yellow squares, $a=0.001$ in blue squares, and $m=4$. The second row: the finite element solution $(u^E_h)_{i,j}$ (left) and the finite difference solution $(u^D_h)_{i,j}$ (right) at all  $(x_i,y_j)$ with $h=1/2^7$. The third row: $\big|(u^E_{h})_{i,j}-(u^E_{h/2})_{2i,2j}\big|$ (left) and  $\big|(u^D_{h})_{i,j}-(u^D_{h/2})_{2i,2j}\big|$ (right) at all  $(x_i,y_j)$ with $h=1/2^6$.}
\label{fig:exam3}
\end{figure}

\begin{figure}[htbp]
	\centering
	\begin{subfigure}[b]{0.4\textwidth}
	\hspace{-1cm}
	\includegraphics[width=7cm,height=7cm]{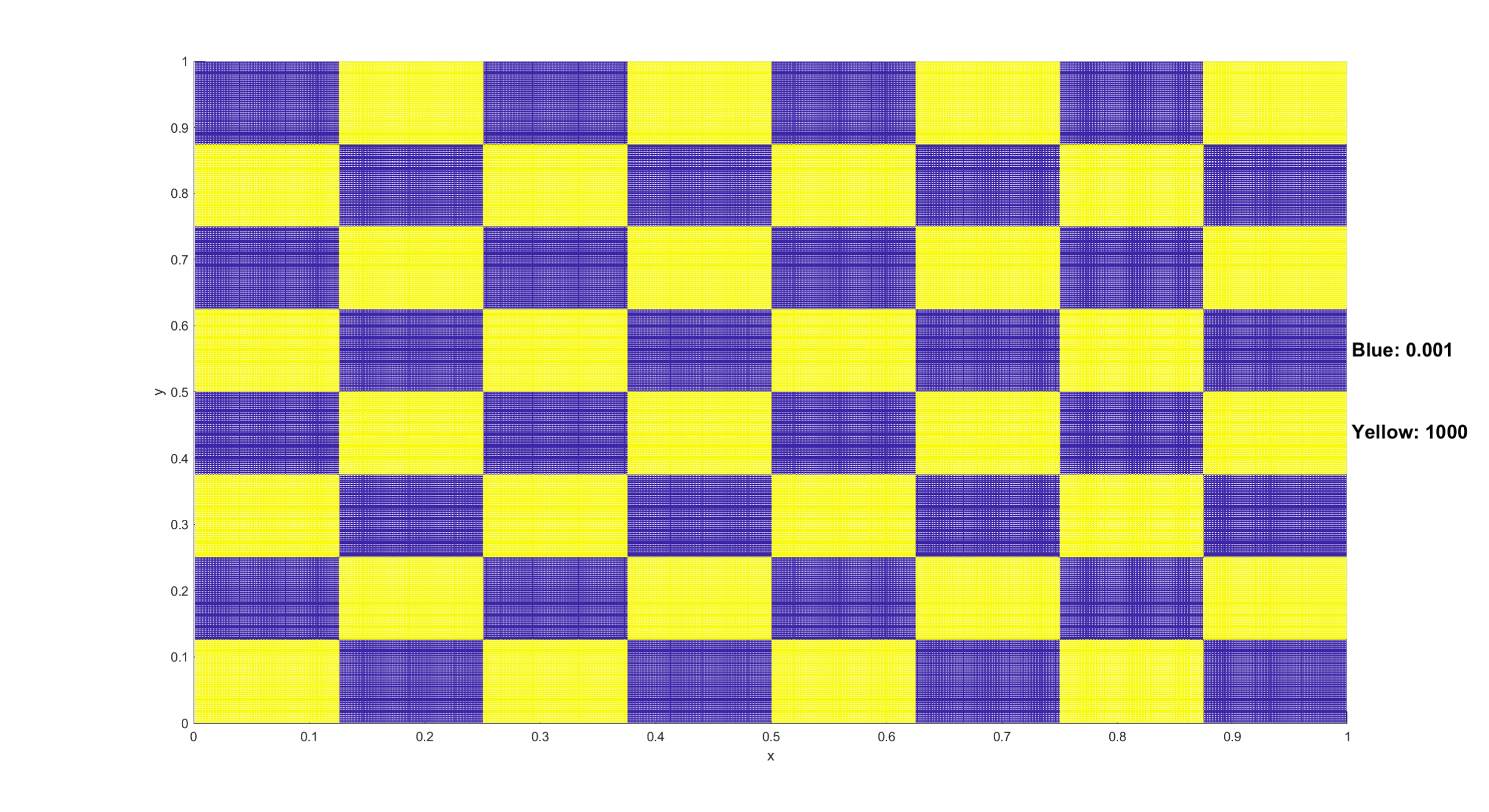}
\end{subfigure}
\begin{subfigure}[b]{0.4\textwidth}
	\hspace{0.7cm}	
	\vspace{0.8cm}	
	\includegraphics[width=5.6cm,height=5.6cm]{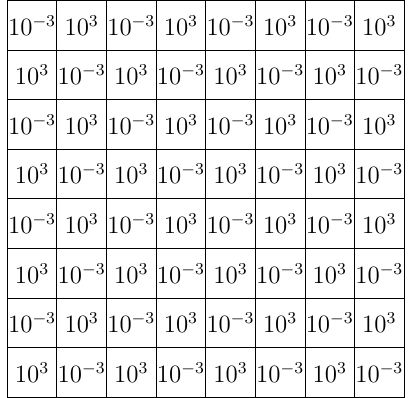}
\end{subfigure}
	\begin{subfigure}[b]{0.45\textwidth}
		\includegraphics[width=7cm,height=7cm]{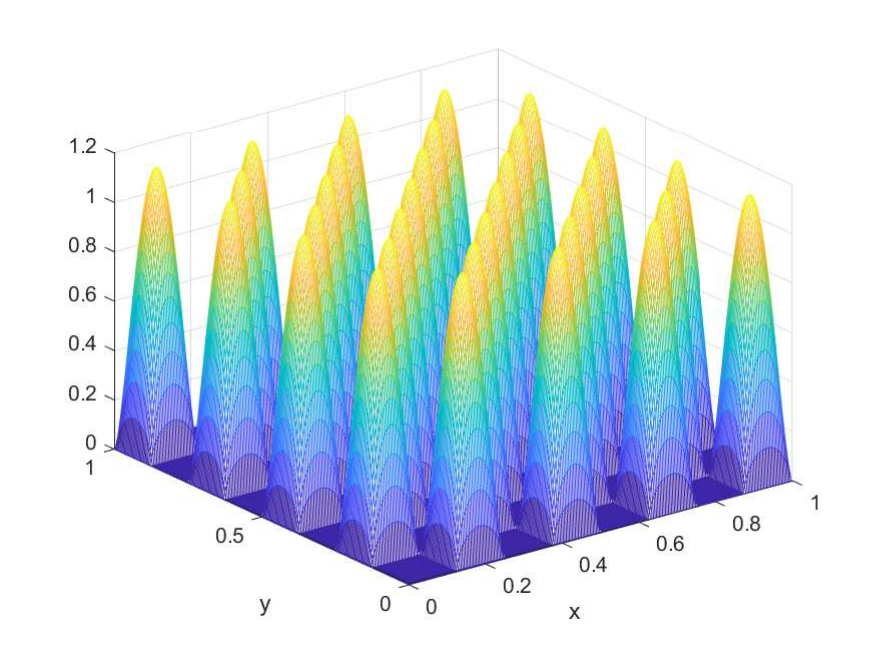}
	\end{subfigure}
	\begin{subfigure}[b]{0.45\textwidth}
		\includegraphics[width=7cm,height=7cm]{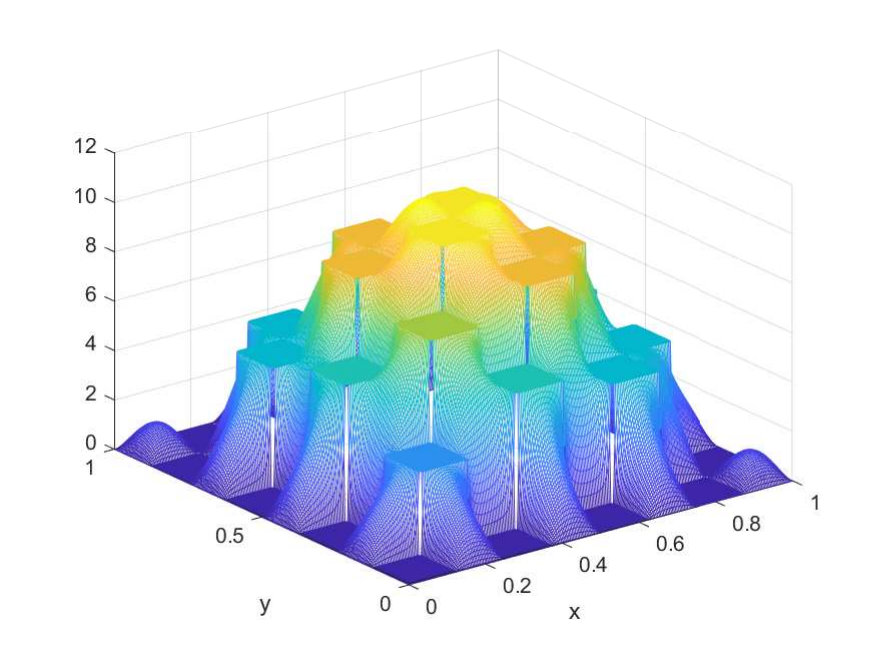}
	\end{subfigure}
	\begin{subfigure}[b]{0.45\textwidth}
		\includegraphics[width=7cm,height=7cm]{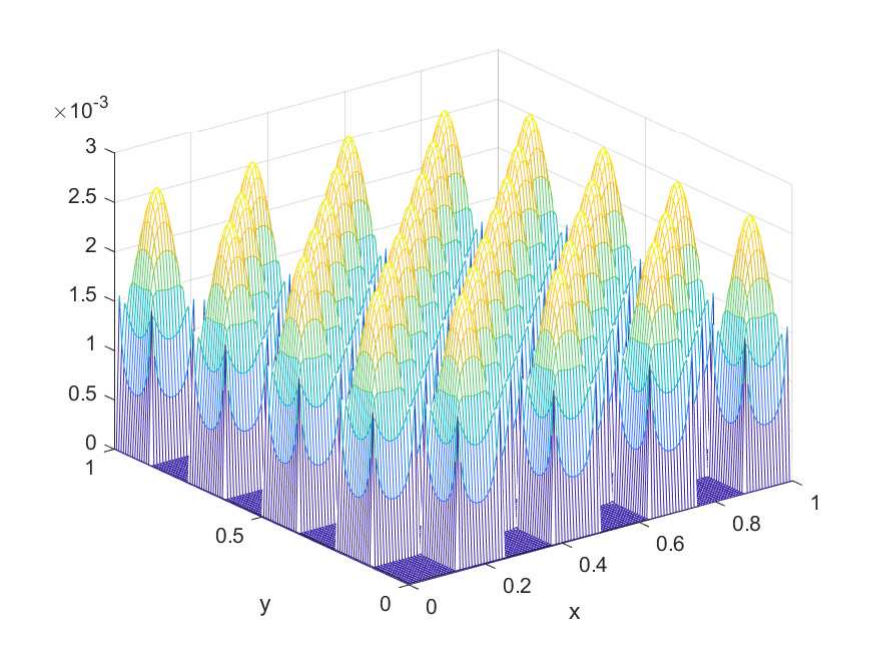}
	\end{subfigure}
	\begin{subfigure}[b]{0.45\textwidth}
		\includegraphics[width=7cm,height=7cm]{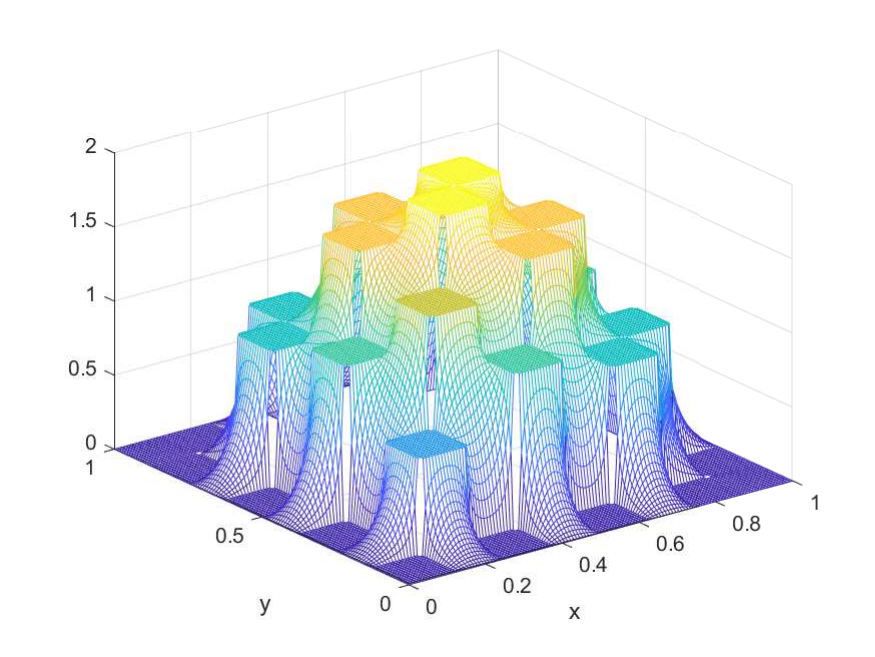}
	\end{subfigure}
	\caption{ \cref{Intersect:ex4}: The first row: the coefficient function $a$ with $a=1000$ in yellow squares, $a=0.001$ in blue squares, and $m=8$. The second row: the finite element solution $(u^E_h)_{i,j}$ (left) and the finite difference solution $(u^D_h)_{i,j}$ (right) at all  $(x_i,y_j)$ with $h=1/2^8$. The third row: $\big|(u^E_{h})_{i,j}-(u^E_{h/2})_{2i,2j}\big|$ (left) and  $\big|(u^D_{h})_{i,j}-(u^D_{h/2})_{2i,2j}\big|$ (right) at all  $(x_i,y_j)$ with $h=1/2^7$.}
\label{fig:exam4}
\end{figure}

\begin{figure}[htbp]
	\centering
	\begin{subfigure}[b]{0.8\textwidth}
		\hspace{3cm}
		\includegraphics[width=7cm,height=7cm]{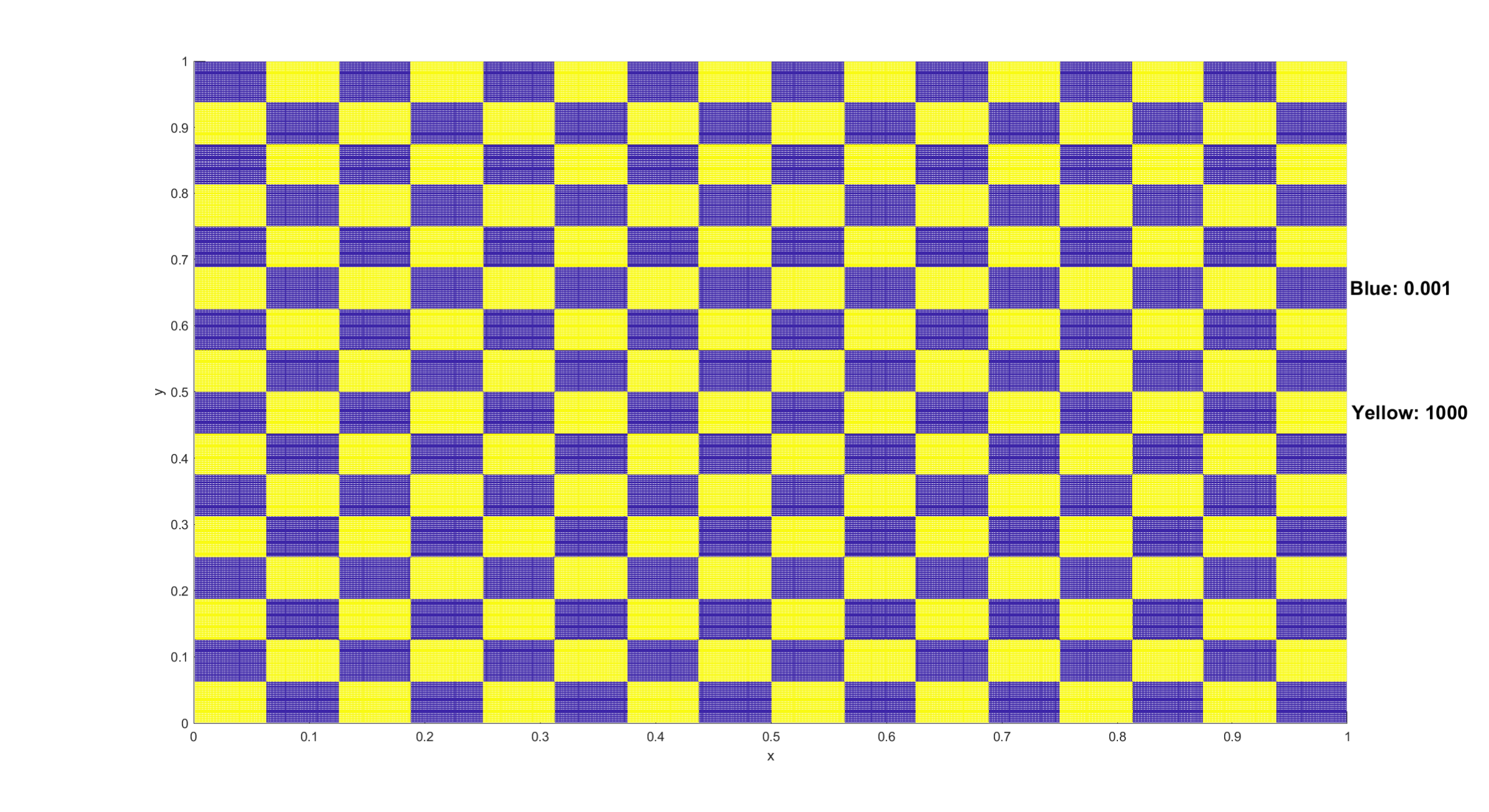}
	\end{subfigure}
	\begin{subfigure}[b]{0.45\textwidth}
		\includegraphics[width=7cm,height=7cm]{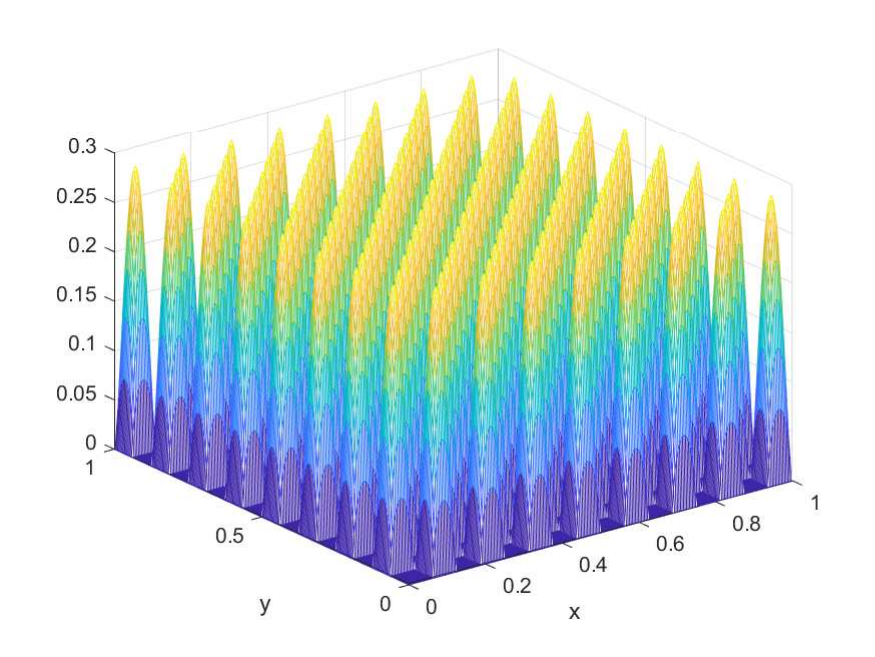}
	\end{subfigure}
	\begin{subfigure}[b]{0.45\textwidth}
		\includegraphics[width=7cm,height=7cm]{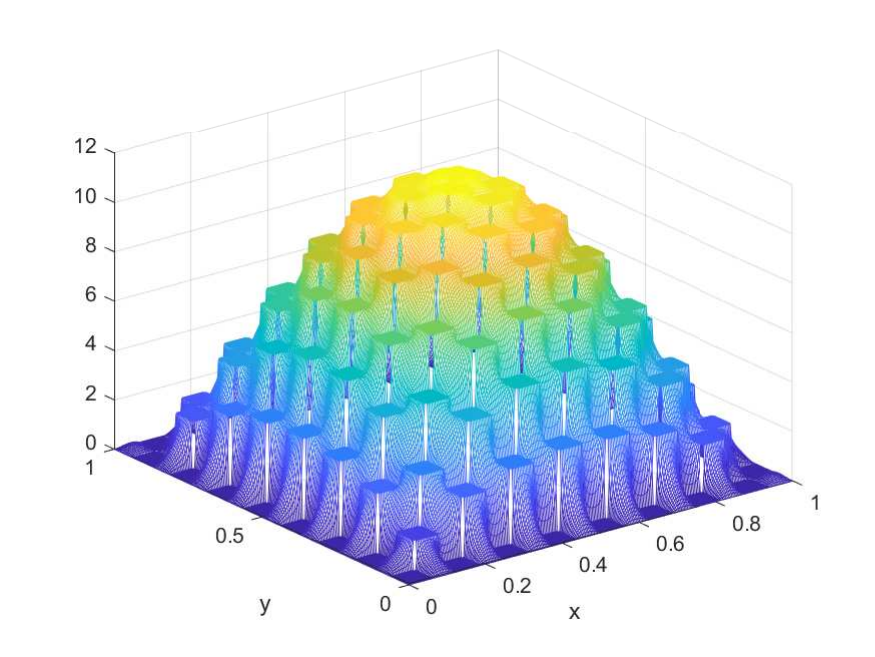}
	\end{subfigure}
	\begin{subfigure}[b]{0.45\textwidth}
		\includegraphics[width=7cm,height=7cm]{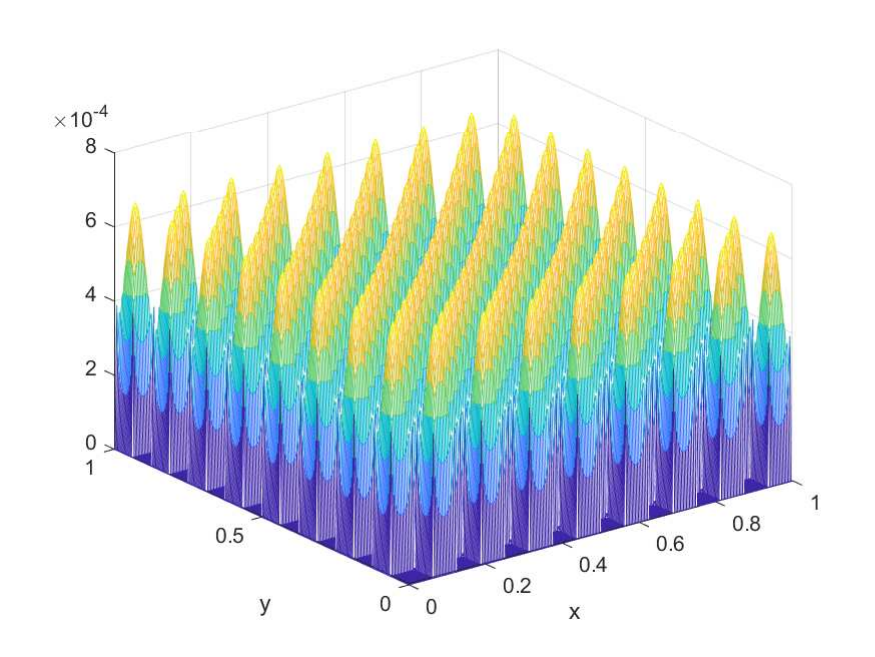}
	\end{subfigure}
	\begin{subfigure}[b]{0.45\textwidth}
		\includegraphics[width=7cm,height=7cm]{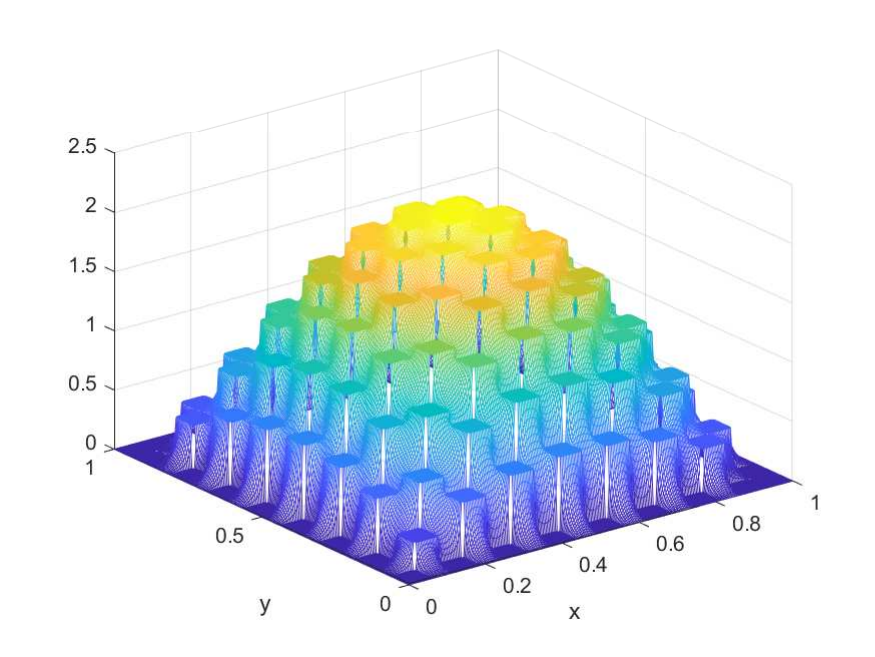}
	\end{subfigure}
	\caption{ \cref{Intersect:ex5}: The first row: the coefficient function $a$ with $a=1000$ in yellow squares, $a=0.001$ in blue squares, and $m=16$. The second row: the finite element solution $(u^E_h)_{i,j}$ (left) and the finite difference solution $(u^D_h)_{i,j}$ (right) at all  $(x_i,y_j)$ with $h=1/2^8$. The third row: $\big|(u^E_{h})_{i,j}-(u^E_{h/2})_{2i,2j}\big|$ (left) and  $\big|(u^D_{h})_{i,j}-(u^D_{h/2})_{2i,2j}\big|$ (right) at all  $(x_i,y_j)$ with $h=1/2^8$.}
\label{fig:exam5}
\end{figure}

\begin{figure}[htbp]
	\centering
	\begin{subfigure}[b]{0.4\textwidth}
	\hspace{-1cm}
	\includegraphics[width=7cm,height=7cm]{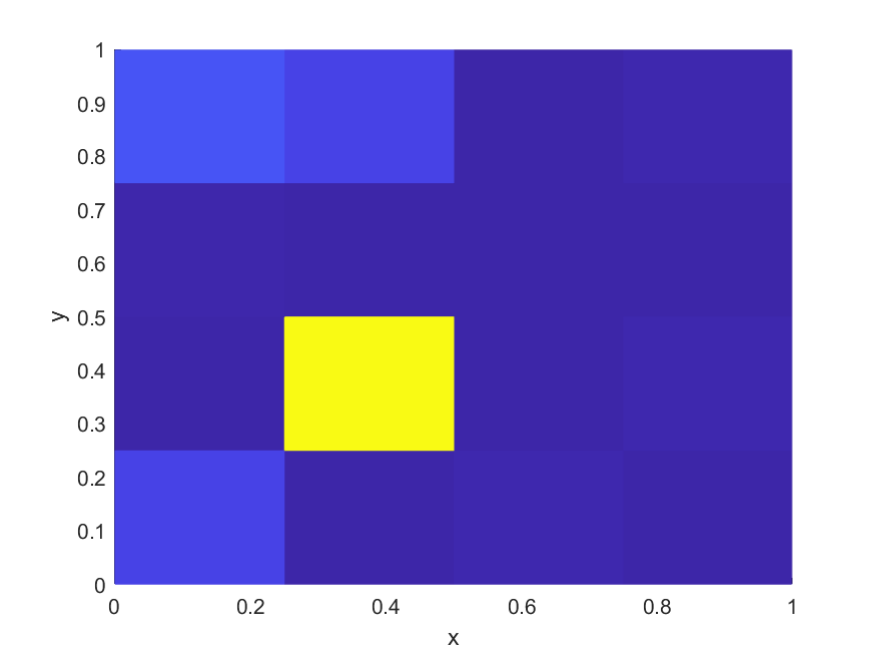}
\end{subfigure}
\begin{subfigure}[b]{0.4\textwidth}
	\hspace{0.7cm}	
	\vspace{0.8cm}	
	\includegraphics[width=5.6cm,height=5.6cm]{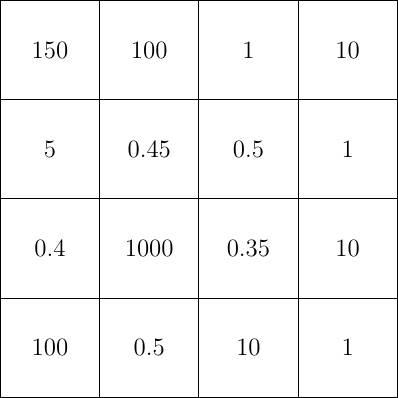}
\end{subfigure}
	\begin{subfigure}[b]{0.45\textwidth}
		\includegraphics[width=7cm,height=7cm]{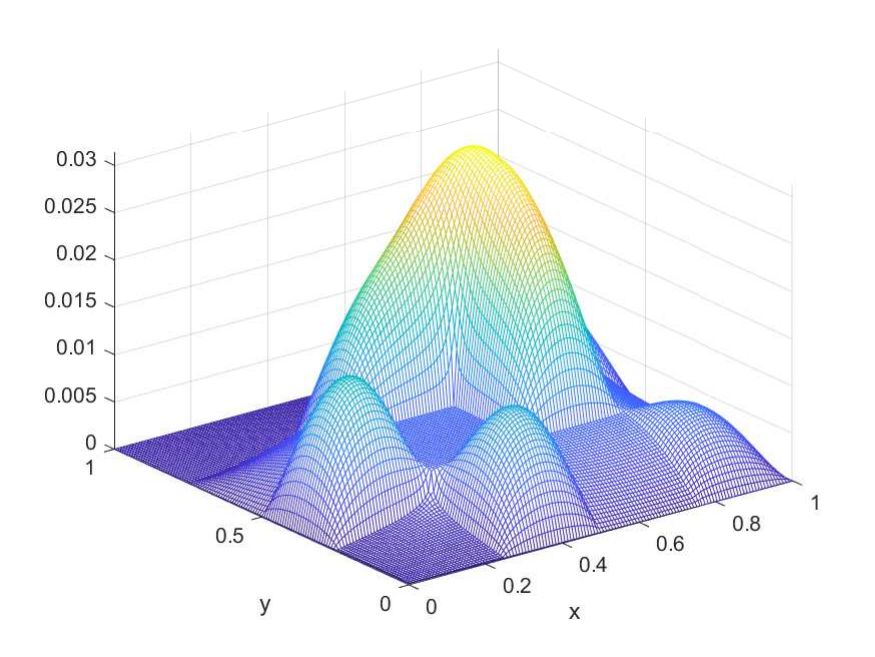}
	\end{subfigure}
	\begin{subfigure}[b]{0.45\textwidth}
		\includegraphics[width=7cm,height=7cm]{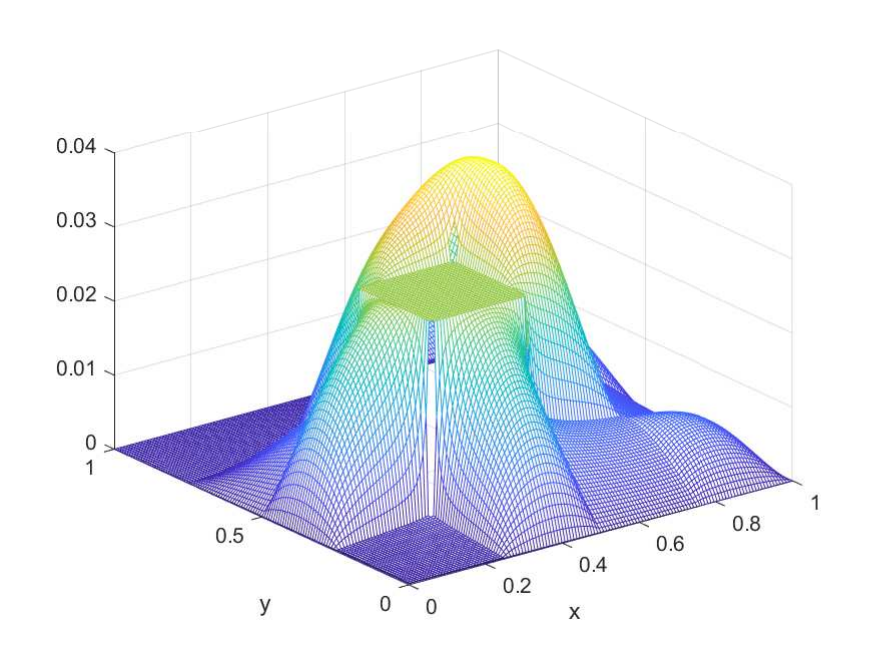}
	\end{subfigure}
	\begin{subfigure}[b]{0.45\textwidth}
		\includegraphics[width=7cm,height=7cm]{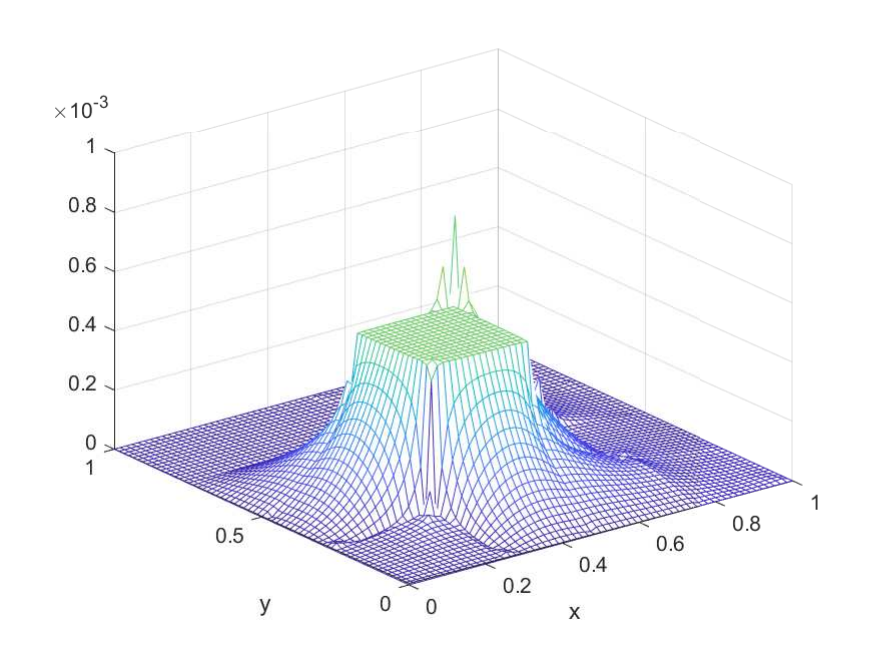}
	\end{subfigure}
	\begin{subfigure}[b]{0.45\textwidth}
		\includegraphics[width=7cm,height=7cm]{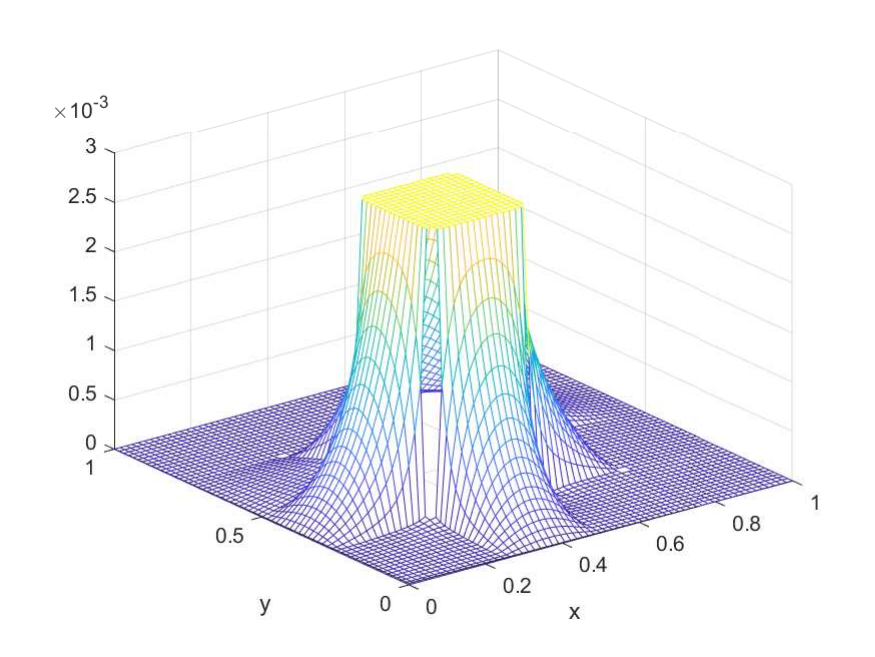}
	\end{subfigure}
	\caption{ \cref{Intersect:ex7}: The first row: the coefficient function $a$ with $m=4$. The second row: the finite element solution $(u^E_h)_{i,j}$ (left) and the finite difference solution $(u^D_h)_{i,j}$ (right) at all  $(x_i,y_j)$ with $h=1/2^7$. The third row: $\big|(u^E_{h})_{i,j}-(u^E_{h/2})_{2i,2j}\big|$ (left) and  $\big|(u^D_{h})_{i,j}-(u^D_{h/2})_{2i,2j}\big|$ (right) at all  $(x_i,y_j)$ with $h=1/2^6$.}
	\label{fig:exam7}
\end{figure}

We shall evaluate finite element and finite difference methods in the $l_2$  norm of the errors given by:
\begin{align*}
 \left\|u^E_{h}-u^E_{h/2}\right\|_{2}:= h\sqrt{\sum_{i=0}^{N}\sum_{j=0}^{N} \left|\left(u^E_{h}\right)_{i,j}-\left(u^E_{h/2}\right)_{2i,2j}\right|^2},
\end{align*}
and
\begin{align*}
	\left\|u^D_{h}-u^D_{h/2}\right\|_{2}:= h\sqrt{{\sum_{i=0 \atop }^{N}\sum_{j=0 \atop }^{N}}_{\hspace{-1.3cm} (x_i,y_j)\notin  \overline{\Gamma^x}\cap \overline{\Gamma^y}} \left|\left(u^D_{h}\right)_{i,j}-\left(u^D_{h/2}\right)_{2i,2j}\right|^2}.
\end{align*}
Furthermore, we also provide results for the infinity norm of the errors given by:
\[
 \left\|u^E_{h}-u^E_{h/2}\right\|_\infty:=\max_{0\le i\le N,\ 0\le j\le N} \left|\left(u^E_{h}\right)_{i,j}-\left(u^E_{h/2}\right)_{2i,2j}\right|,
\]
and
\[
\left\|u^D_{h}-u^D_{h/2}\right\|_\infty:=\max_{0\le i\le N,\ 0\le j\le N \atop   (x_i,y_j)\notin  \overline{\Gamma^x}\cap \overline{\Gamma^y} } \left|\left(u^D_{h}\right)_{i,j}-\left(u^D_{h/2}\right)_{2i,2j}\right|.
\]
%
%
%
%
\begin{example}\label{Intersect:ex1}
	\normalfont
	Let $f=1$ and $m=2$ in \eqref{intersect:3}.
	The coefficient function $a$  in \eqref{intersect:3} is given in the first row of \cref{fig:exam1}.  Note that $a=1000$ in yellow squares and $a=0.001$ in blue squares.
	The numerical results are displayed in \cref{Intersect:table1} and \cref{fig:exam1}.	
\end{example}
\begin{table}[htbp]
	\caption{Performance in \cref{Intersect:ex1} of finite element and finite difference methods on uniform Cartesian meshes. }
	\centering
 \renewcommand{\arraystretch}{1.5}
	\setlength{\tabcolsep}{0.3mm}{
		\begin{tabular}{c|c|c|c|c|c|c|c|c}
			\hline
			\multicolumn{1}{c|}{}  &
			\multicolumn{4}{c|}{Use FEM in \cref{Sec:FEM}}  &
            \multicolumn{4}{c}{Use FDM in \cref{Sec:FDM}} \\
			\hline
			$h$
			&   $\big\|u^E_{h}-u^E_{h/2}\big\|_{2}$
			&order &  $\big\|u^E_{h}-u^E_{h/2}\big\|_\infty$ & order &   $\big\|u^D_{h}-u^D_{h/2}\big\|_{2}$
			&order &  $\big\|u^D_{h}-u^D_{h/2}\big\|_\infty$ & order \\
			\hline
$\frac{1}{2^2}$  &1.4205E+00   &   &4.0179E+00   &   &6.9053E-01   &   &1.9531E+00   & \\
$\frac{1}{2^3}$  &3.1139E-01   &2.19   &7.7007E-01   &2.38   &2.5329E-01   &1.45   &6.1753E-01   &1.66\\
$\frac{1}{2^4}$  &7.3853E-02   &2.08   &1.7475E-01   &2.14   &6.9442E-02   &1.87   &1.6578E-01   &1.90\\
$\frac{1}{2^5}$  &1.8132E-02   &2.03   &4.2797E-02   &2.03   &1.7806E-02   &1.96   &4.2243E-02   &1.97\\
$\frac{1}{2^6}$  &4.5073E-03   &2.01   &1.0647E-02   &2.01   &4.4838E-03   &1.99   &1.0612E-02   &1.99\\			
			\hline
	\end{tabular}}
	\label{Intersect:table1}
\end{table}	
\begin{example}\label{Intersect:ex2}
	\normalfont
	Let $f=1$ and $m=4$ in \eqref{intersect:3}.
	The coefficient function $a$  in \eqref{intersect:3} is given in the first row of \cref{fig:exam2}. Note that $a=10$ in yellow squares and $a=1$ in blue squares.
	The numerical results are displayed in \cref{Intersect:table2} and \cref{fig:exam2}.	
\end{example}
\begin{table}[htbp]
	\caption{Performance in \cref{Intersect:ex2} of finite element and finite difference methods on uniform Cartesian meshes. }
	\centering
	 \renewcommand{\arraystretch}{1.5}
	\setlength{\tabcolsep}{0.3mm}{
		\begin{tabular}{c|c|c|c|c|c|c|c|c}
			\hline
			\multicolumn{1}{c|}{}  &
			\multicolumn{4}{c|}{Use FEM in \cref{Sec:FEM}}  &
\multicolumn{4}{c}{Use FDM in \cref{Sec:FDM}} \\
			\hline
			$h$
			&   $\big\|u^E_{h}-u^E_{h/2}\big\|_{2}$
			&order &  $\big\|u^E_{h}-u^E_{h/2}\big\|_\infty$ & order &   $\big\|u^D_{h}-u^D_{h/2}\big\|_{2}$
			&order &  $\big\|u^D_{h}-u^D_{h/2}\big\|_\infty$ & order \\
			\hline
$\frac{1}{2^3}$   &5.0839E-04   &   &1.5140E-03   &   &1.9897E-03   &   &4.9999E-03   & \\
$\frac{1}{2^4}$   &3.1692E-04   &0.68   &8.6841E-04   &0.80   &1.0171E-03   &0.97   &2.4015E-03   &1.06\\
$\frac{1}{2^5}$   &2.0231E-04   &0.65   &5.8430E-04   &0.57   &5.4431E-04   &0.90   &1.3437E-03   &0.84\\
$\frac{1}{2^6}$   &1.2417E-04   &0.70   &3.9910E-04   &0.55   &3.0031E-04   &0.86   &8.2004E-04   &0.71\\		
			\hline
	\end{tabular}}
	\label{Intersect:table2}
\end{table}	
\begin{example}\label{Intersect:ex3:New}
	\normalfont
	Let $f=1$ and $m=8$ in \eqref{intersect:3}.
	The coefficient function $a$  in \eqref{intersect:3} is given in the first row of \cref{fig:exam3:New}.  Note that $a=10$ in yellow squares and $a=1$ in blue squares.
	The numerical results are displayed in \cref{Intersect:table3:New} and \cref{fig:exam3:New}.	
\end{example}
\begin{table}[htbp]
	\caption{Performance in \cref{Intersect:ex3:New} of finite element and finite difference methods on uniform Cartesian meshes. }
	\centering
	\renewcommand{\arraystretch}{1.5}
	\setlength{\tabcolsep}{0.3mm}{
		\begin{tabular}{c|c|c|c|c|c|c|c|c}
			\hline
			\multicolumn{1}{c|}{}  &
			\multicolumn{4}{c|}{Use FEM in \cref{Sec:FEM}}  &
			\multicolumn{4}{c}{Use FDM in \cref{Sec:FDM}} \\
			\hline
			$h$
			&   $\big\|u^E_{h}-u^E_{h/2}\big\|_{2}$
			&order &  $\big\|u^E_{h}-u^E_{h/2}\big\|_\infty$ & order &   $\big\|u^D_{h}-u^D_{h/2}\big\|_{2}$
			&order &  $\big\|u^D_{h}-u^D_{h/2}\big\|_\infty$ & order \\
			\hline
$\frac{1}{2^4}$   &8.0188E-04   &   &1.7532E-03   &   &2.8484E-03   &   &6.7857E-03   & \\
$\frac{1}{2^5}$   &4.9292E-04   &0.70   &1.0323E-03   &0.76   &1.4258E-03   &1.00   &2.9368E-03   &1.21\\
$\frac{1}{2^6}$   &3.0106E-04   &0.71   &6.3921E-04   &0.69   &7.7563E-04   &0.88   &1.5856E-03   &0.89\\
$\frac{1}{2^7}$   &1.8072E-04   &0.74   &3.9633E-04   &0.69   &4.2936E-04   &0.85   &8.9909E-04   &0.82\\	
			\hline
	\end{tabular}}
	\label{Intersect:table3:New}
\end{table}	
\begin{example}\label{Intersect:ex3}
	\normalfont
	Let $f=1$ and $m=4$ in \eqref{intersect:3}.
	The coefficient function $a$  in \eqref{intersect:3} is given in the first row of \cref{fig:exam3}.  Note that $a=1000$ in yellow squares and $a=0.001$ in blue squares.
	The numerical results are displayed in \cref{Intersect:table3} and \cref{fig:exam3}.	
\end{example}
\begin{table}[htbp]
	\caption{Performance in \cref{Intersect:ex3} of finite element and finite difference methods on uniform Cartesian meshes. }
	\centering
	 \renewcommand{\arraystretch}{1.5}
	\setlength{\tabcolsep}{0.3mm}{
		\begin{tabular}{c|c|c|c|c|c|c|c|c}
			\hline
			\multicolumn{1}{c|}{}  &
			\multicolumn{4}{c|}{Use FEM in \cref{Sec:FEM}}  &
\multicolumn{4}{c}{Use FDM in \cref{Sec:FDM}} \\
			\hline
			$h$
			&   $\big\|u^E_{h}-u^E_{h/2}\big\|_{2}$
			&order &  $\big\|u^E_{h}-u^E_{h/2}\big\|_\infty$ & order &   $\big\|u^D_{h}-u^D_{h/2}\big\|_{2}$
			&order &  $\big\|u^D_{h}-u^D_{h/2}\big\|_\infty$ & order \\
			\hline
$\frac{1}{2^3}$   &3.5513E-01   &   &1.0045E+00   &   &2.5237E+00   &   &6.2348E+00   & \\
$\frac{1}{2^4}$   &7.7843E-02   &2.19   &1.9252E-01   &2.38   &1.5017E+00   &0.75   &3.4263E+00   &0.86\\
$\frac{1}{2^5}$   &1.8458E-02   &2.08   &4.3687E-02   &2.14   &9.5379E-01   &0.65   &2.1370E+00   &0.68\\
$\frac{1}{2^6}$   &4.5272E-03   &2.03   &1.0699E-02   &2.03   &6.5185E-01   &0.55   &1.4513E+00   &0.56\\			
			\hline
	\end{tabular}}
	\label{Intersect:table3}
\end{table}	
\begin{example}\label{Intersect:ex4}
	\normalfont
	Let $f=1$ and $m=8$ in \eqref{intersect:3}.
	The coefficient function $a$  in \eqref{intersect:3} is given in the first row of \cref{fig:exam4}.  Note that $a=1000$ in yellow squares and $a=0.001$ in blue squares.
	The numerical results are displayed in \cref{Intersect:table4} and \cref{fig:exam4}.	
\end{example}
\begin{table}[htbp]
	\caption{Performance in \cref{Intersect:ex4} of finite element and finite difference methods on uniform Cartesian meshes. }
	\centering
	 \renewcommand{\arraystretch}{1.5}
	\setlength{\tabcolsep}{0.3mm}{
		\begin{tabular}{c|c|c|c|c|c|c|c|c}
			\hline
			\multicolumn{1}{c|}{}  &
			\multicolumn{4}{c|}{Use FEM in \cref{Sec:FEM}}  &
\multicolumn{4}{c}{Use FDM in \cref{Sec:FDM}} \\
			\hline
			$h$
			&   $\big\|u^E_{h}-u^E_{h/2}\big\|_{2}$
			&order &  $\big\|u^E_{h}-u^E_{h/2}\big\|_\infty$ & order &   $\big\|u^D_{h}-u^D_{h/2}\big\|_{2}$
			&order &  $\big\|u^D_{h}-u^D_{h/2}\big\|_\infty$ & order \\
			\hline
$\frac{1}{2^4}$   &8.8776E-02   &   &2.5112E-01   &   &3.8360E+00   &   &9.1066E+00   & \\
$\frac{1}{2^5}$   &1.9452E-02   &2.19   &4.8129E-02   &2.38   &2.1341E+00   &0.85   &4.3654E+00   &1.06\\
$\frac{1}{2^6}$   &4.6051E-03   &2.08   &1.0922E-02   &2.14   &1.3691E+00   &0.64   &2.7437E+00   &0.67\\
$\frac{1}{2^7}$   &1.1223E-03   &2.04   &2.6748E-03   &2.03   &9.4162E-01   &0.54   &1.8822E+00   &0.54\\			
			\hline
	\end{tabular}}
	\label{Intersect:table4}
\end{table}	
\begin{example}\label{Intersect:ex5}
	\normalfont
	Let $f=1$ and $m=16$ in \eqref{intersect:3}.
	The coefficient function $a$  in \eqref{intersect:3} is given in the first row of \cref{fig:exam5}.  Note that $a=1000$ in yellow squares and $a=0.001$ in blue squares.
	The numerical results are displayed in \cref{Intersect:table5} and \cref{fig:exam5}.	
\end{example}
\begin{table}[htbp]
	\caption{Performance in \cref{Intersect:ex5} of finite element and finite difference methods on uniform Cartesian meshes. }
	\centering
	 \renewcommand{\arraystretch}{1.5}
	\setlength{\tabcolsep}{0.3mm}{
		\begin{tabular}{c|c|c|c|c|c|c|c|c}
			\hline
			\multicolumn{1}{c|}{}  &
			\multicolumn{4}{c|}{Use FEM in \cref{Sec:FEM}}  &
\multicolumn{4}{c}{Use FDM in \cref{Sec:FDM}} \\
			\hline
			$h$
			&   $\big\|u^E_{h}-u^E_{h/2}\big\|_{2}$
			&order &  $\big\|u^E_{h}-u^E_{h/2}\big\|_\infty$ & order &   $\big\|u^D_{h}-u^D_{h/2}\big\|_{2}$
			&order &  $\big\|u^D_{h}-u^D_{h/2}\big\|_\infty$ & order \\
			\hline
$\frac{1}{2^5}$   &2.2187E-02   &   &6.2779E-02   &   &4.7871E+00   &   &1.0605E+01   & \\
$\frac{1}{2^6}$   &4.8528E-03   &2.19   &1.2032E-02   &2.38   &2.5505E+00   &0.91   &4.9618E+00   &1.10\\
$\frac{1}{2^7}$   &1.1399E-03   &2.09   &2.7305E-03   &2.14   &1.6425E+00   &0.63   &3.1229E+00   &0.67\\
$\frac{1}{2^8}$   &2.6975E-04   &2.08   &6.6872E-04   &2.03   &1.1339E+00   &0.53   &2.1465E+00   &0.54\\			
			\hline
	\end{tabular}}
	\label{Intersect:table5}
\end{table}	
\begin{example}\label{Intersect:ex7}
	\normalfont
	Let $f=1$ and $m=4$ in \eqref{intersect:3}.
	The coefficient function $a$  in \eqref{intersect:3} is given in the first row of \cref{fig:exam7}.  
	The numerical results are displayed in \cref{Intersect:table7} and \cref{fig:exam7}.	
\end{example}
\begin{table}[htbp]
	\caption{Performance in \cref{Intersect:ex7} of finite element and finite difference methods on uniform Cartesian meshes. }
	\centering
	\renewcommand{\arraystretch}{1.5}
	\setlength{\tabcolsep}{0.3mm}{
		\begin{tabular}{c|c|c|c|c|c|c|c|c}
			\hline
			\multicolumn{1}{c|}{}  &
			\multicolumn{4}{c|}{Use FEM in \cref{Sec:FEM}}  &
			\multicolumn{4}{c}{Use FDM in \cref{Sec:FDM}} \\
			\hline
			$h$
			&   $\big\|u^E_{h}-u^E_{h/2}\big\|_{2}$
			&order &  $\big\|u^E_{h}-u^E_{h/2}\big\|_\infty$ & order &   $\big\|u^D_{h}-u^D_{h/2}\big\|_{2}$
			&order &  $\big\|u^D_{h}-u^D_{h/2}\big\|_\infty$ & order \\
			\hline
$\frac{1}{2^3}$   &5.7156E-04   &   &2.8756E-03   &   &3.9260E-03   &   &1.1793E-02   & \\
$\frac{1}{2^4}$   &2.5777E-04   &1.15   &1.3236E-03   &1.12   &2.2782E-03   &0.79   &6.7839E-03   &0.80\\
$\frac{1}{2^5}$   &2.1484E-04   &0.26   &1.0128E-03   &0.39   &1.4528E-03   &0.65   &4.3152E-03   &0.65\\
$\frac{1}{2^6}$   &1.9528E-04   &0.14   &8.2897E-04   &0.29   &9.9857E-04   &0.54   &2.9800E-03   &0.53\\		
			\hline
	\end{tabular}}
	\label{Intersect:table7}
\end{table}	
\section{Possible explanations}\label{Intersect:sec:Reason}

\begin{figure}[htbp]
	\centering
	\begin{subfigure}[b]{0.45\textwidth}
		\includegraphics[width=7.5cm,height=7.5cm]{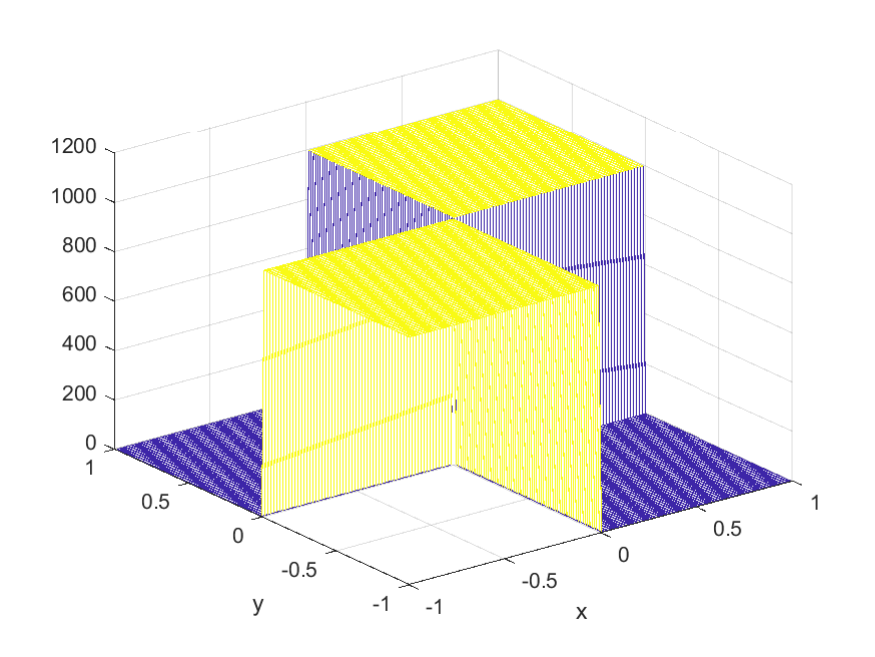}
	\end{subfigure}
	\begin{subfigure}[b]{0.45\textwidth}
		\includegraphics[width=7.5cm,height=7.5cm]{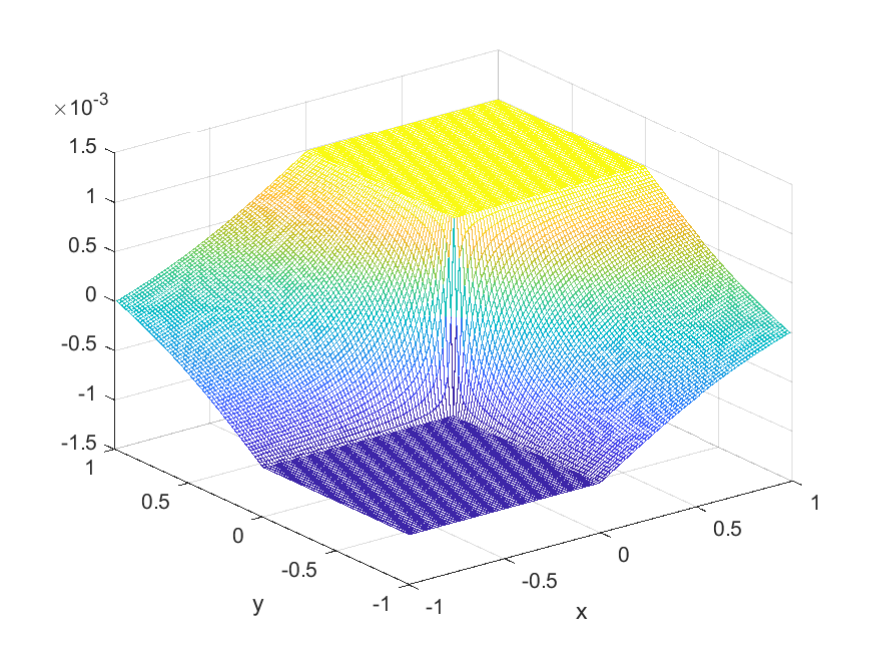}
	\end{subfigure}
	\begin{subfigure}[b]{0.45\textwidth}
		\includegraphics[width=7.5cm,height=7.5cm]{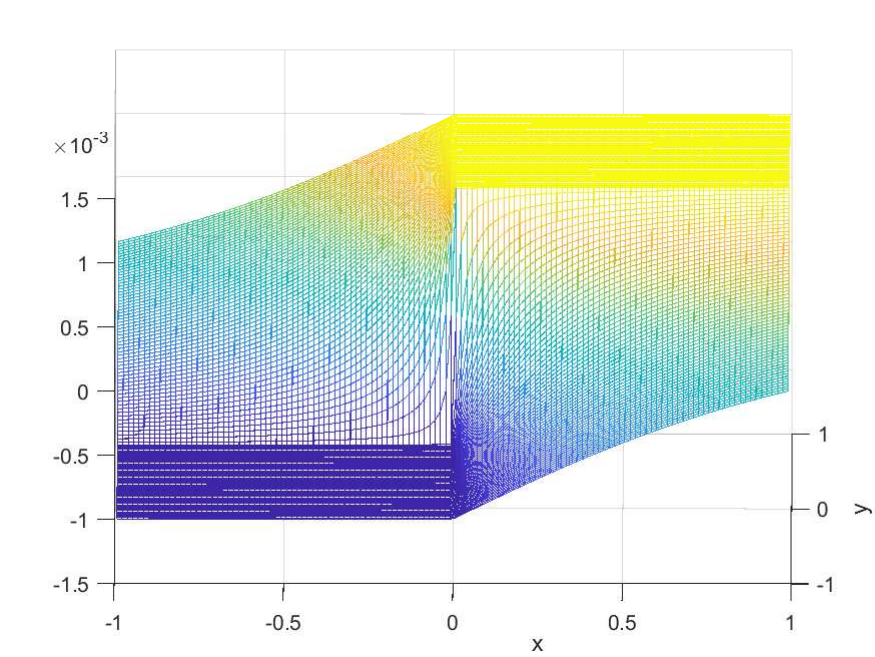}
	\end{subfigure}
	\begin{subfigure}[b]{0.45\textwidth}
		\includegraphics[width=7.5cm,height=7.5cm]{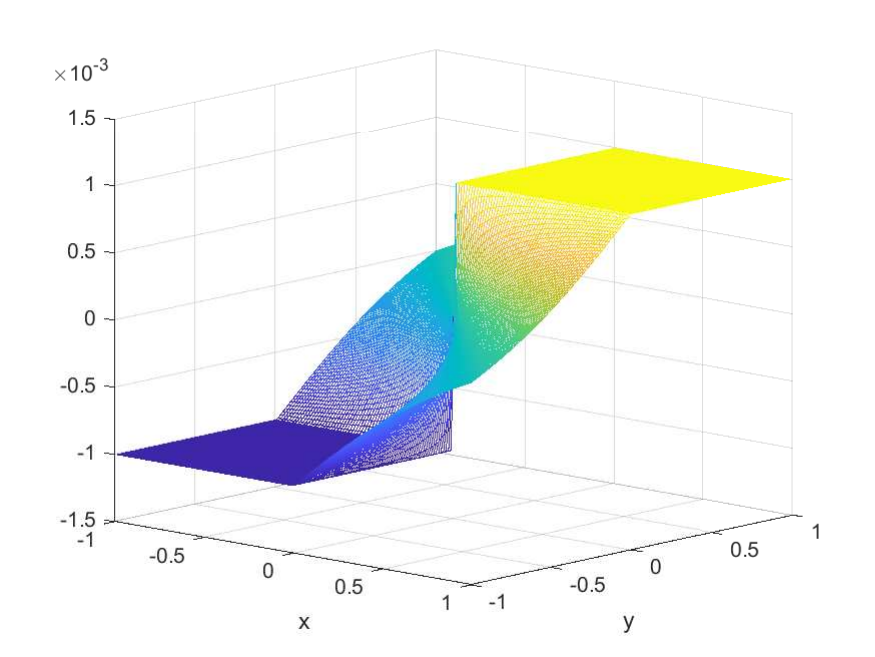}
	\end{subfigure}
	\caption{The first row: the coefficient function $a$ with $a=1000$ in yellow squares, $a=0.001$ in blue squares (left) and the analytical solution $u$ in \eqref{exact:u} with $\lambda$ in \eqref{Exact_lam} (right). The second row:  the analytical solution $u$ in \eqref{exact:u} with $\lambda$ in \eqref{Exact_lam}.}
	\label{fig:ExactU}
\end{figure}

The reason why the finite element method and finite difference method generate two distinct numerical solutions might be due to the severe singularities caused by high jumps in the coefficient functions at the intersection points (see e.g. \cite{CaKi01,Kell75,KiKo06,Blumen85,Petz01,Petz02,Nic90,KCPK07,NS94,FGW73,Kell71,Kell72}).
Recall that $\cup_{p=1}^m\cup_{q=1}^m \Omega_{p,q}=\Omega\setminus \Gamma$ and $a$ is a constant in each $\Omega_{p,q}$. Using \cite[Assumption 4.43]{PietroErn2012}, we define that $P_{\Omega}:=\{\Omega_{p,q}\}_{1\le p,q\le m}$, i.e., $P_{\Omega}$ is the partition of $\Omega$.
By the discussions in  \cite{CaKi01,Kell75,KiKo06,Blumen85,Petz01,Petz02,Nic90,KCPK07,NS94,FGW73,Kell71,Kell72}, we have that the exact solution $u$ of the model problem \eqref{intersect:3} may satisfy  $u\in H^{1+\lambda}(P_{\Omega})$ with $0<\lambda<1$.
 For example: Let $\Omega=(-1,1)^2$, $0<\lambda<1$,  $\theta=\arctan\left(\frac{y}{x}\right)$,
 \be\label{Exact:a}
 a=
 \begin{cases}
 	\ka^{-1}& 0< \theta< \pi/2,\\
 	\ka&\pi/2< \theta< \pi,\\
 	\ka^{-1}&\pi< \theta< 3\pi/2,\\
 	\ka&3\pi/2< \theta< 2\pi,
 \end{cases}
 \qquad\qquad \ka=\tan(\lambda\pi/4),  
 \ee
 \[
 \begin{split}
 	 \Gamma:=\Big\{(x,y) : x=0,\ y\in (-1,0)\cup(0,1)\Big \} \cup \Big\{(x,y) : y=0,\ x\in (-1,0)\cup(0,1)\Big \},
 \end{split}
 \]
and consider the following problem with a zero source term and the non-homogeneous  Dirichlet boundary condition:
 \be\label{problem:new}
 \left\{ \begin{array}{llll}
 	-\nabla \cdot \Big(a \nabla  u\Big)&=&0 
 	&\mbox{in } \Omega \setminus {\Gamma},\\
 	\left[u\right]&=&0 
 	&\mbox{on }\Gamma,\\	
 	\left[a \nabla  u \cdot \vec{n}\right] &=&0 
 	&\mbox{on }\Gamma,\\
 	u&=&g 
 	&\mbox{on }	\partial\Omega.
 \end{array}
 \right.
 \ee
 By \cite{ChenDai2002,Kell75,Petz02}, the analytical solution $u$ of \eqref{problem:new} is (see \cref{fig:ExactU} for an illustration)
\be\label{exact:u}
u=
\begin{cases}
u_1=	(\sqrt{x^2+y^2})^{\lambda}\sin(\lambda \pi/4) \cos\left(\frac{\lambda(\pi-4\theta)}{4} \right)    &\quad 0\le \theta \le \pi/2,\\
u_2=	(\sqrt{x^2+y^2})^{\lambda}\cos(\lambda\pi/4) \sin\left(\frac{\lambda(3\pi-4\theta)}{4} \right)       & \quad\pi/2\le \theta\le \pi,\\
u_3=	-(\sqrt{x^2+y^2})^{\lambda}\sin(\lambda \pi/4) \cos\left(\frac{\lambda(5\pi-4\theta)}{4} \right)  & \quad\pi\le \theta\le 3\pi/2,\\
u_4=-	(\sqrt{x^2+y^2})^{\lambda}\cos(\lambda\pi/4) \sin\left(\frac{\lambda(7\pi-4\theta)}{4} \right)    & \quad3\pi/2\le \theta\le 2\pi, 
\end{cases}
\ee
 where the boundary function $g$ in \eqref{problem:new} can be obtained by plugging $u$ in \eqref{exact:u} into \eqref{problem:new}. Let  $r=\sqrt{x^2+y^2}$. Then
\be\label{exact:ux}
u_x=
\begin{cases}
	(u_1)_x=	\lambda r^{\lambda-1/2}  \sin(\lambda \pi/4) \left\{ \cos\left(\frac{\lambda(\pi-4\theta)}{4} \right) x -\sin\left(\frac{\lambda(\pi-4\theta)}{4} \right)y\right\}     &\quad 0< \theta < \pi/2,\\
(u_2)_x=		\lambda r^{\lambda-1/2}  \cos(\lambda \pi/4) \left\{ \sin\left(\frac{\lambda(3\pi-4\theta)}{4} \right) x + \cos\left(\frac{\lambda(3\pi-4\theta)}{4} \right) y\right\}         &\quad \pi/2< \theta< \pi,\\
	(u_3)_x=	-\lambda r^{\lambda-1/2}  \sin(\lambda \pi/4) \left\{ \cos\left(\frac{\lambda(5\pi-4\theta)}{4} \right) x -\sin\left(\frac{\lambda(5\pi-4\theta)}{4} \right)y\right\}    &\quad \pi< \theta< 3\pi/2,\\
(u_4)_x=-\lambda r^{\lambda-1/2}  \cos(\lambda \pi/4) \left\{ \sin\left(\frac{\lambda(7\pi-4\theta)}{4} \right) x + \cos\left(\frac{\lambda(7\pi-4\theta)}{4} \right) y\right\}    &\quad 3\pi/2< \theta< 2\pi, 
\end{cases}
\ee 
\be\label{exact:uy}
u_y=
\begin{cases}
	(u_1)_y=	\lambda r^{\lambda-1/2}  \sin(\lambda \pi/4) \left\{\sin\left(\frac{\lambda(\pi-4\theta)}{4} \right)x+\cos\left(\frac{\lambda(\pi-4\theta)}{4} \right)y \right\}     &\quad 0< \theta < \pi/2,\\
	(u_2)_y=		\lambda r^{\lambda-1/2}  \cos(\lambda \pi/4) \left\{- \cos\left(\frac{\lambda(3\pi-4\theta)}{4} \right) x + \sin\left(\frac{\lambda(3\pi-4\theta)}{4} \right) y\right\}         &\quad \pi/2< \theta< \pi,\\
	(u_3)_y=-\lambda r^{\lambda-1/2}  \sin(\lambda \pi/4) \left\{\sin\left(\frac{\lambda(5\pi-4\theta)}{4} \right)x+\cos\left(\frac{\lambda(5\pi-4\theta)}{4} \right)y\right\}    &\quad \pi< \theta< 3\pi/2,\\
	(u_4)_y=-	\lambda r^{\lambda-1/2}  \cos(\lambda \pi/4) \left\{- \cos\left(\frac{\lambda(7\pi-4\theta)}{4} \right) x + \sin\left(\frac{\lambda(7\pi-4\theta)}{4} \right) y\right\}   &\quad 3\pi/2< \theta< 2\pi.
\end{cases}
\ee 
Choose
\be\label{Exact_lam}
\lambda=\arctan\left(10^{-3}\right)\frac{4}{\pi}\approx 0.001273239120322.
\ee
We have that $\ka=10^{-3}$ in \eqref{Exact:a}, and the $u_x$ in \eqref{exact:ux} and $u_y$ in \eqref{exact:uy} meet that 
\[
\left|\frac{\partial u_i(0,0)}{\partial x}\right|=\infty, \qquad  \left|\frac{\partial u_i(0,0)}{\partial y}\right|=\infty, \quad \text{for} \quad i=1,2,3,4.
\]
Note that $a$ in \eqref{Exact:a} and $u$ in \eqref{exact:u} with \eqref{Exact_lam} are shown in \cref{fig:ExactU}.

For the model problem \eqref{intersect:3}, \cite[Section 4.5]{PietroErn2012} discussed the error analysis and  convergence proof of  the discontinuous Galerkin method. Specifically, by adding the penalty term and analyzing the regularity, \cite[Theorem 4.53 in Section 4.5.3]{PietroErn2012} proposed the error estimate for the smooth solution and \cite[Theorem 4.56 in Section 4.5.4]{PietroErn2012} provided the error estimate for the low-regularity solution. 
If the exact solution $u$ has uniformly continuous partial derivatives of (total) orders up to seven in each $\Omega_{p,q}$, \cite[Theorem 3.1]{FengHanMinevCross} theoretically proved the sixth-order convergence rate by the discrete maximum principle for the FDMs in  \cite[Theorems 2.1-2.3]{FengHanMinevCross}. However, to the best of our knowledge, no FDMs have been presented in the literature that can theoretically establish the proof of the convergence rate for the model problem \eqref{intersect:3} with $u \in H^{1+\lambda}(P_\Omega)$ and $0<\lambda<1$.

In \cref{Intersect:ex1}, the domain $\Omega$  contains only 4 subdomains, the coefficient function $a$ exhibits high jumps, and the zero Dirichlet boundary condition is imposed. The boundary effect possibly causes  the solution  $u$ to be smooth in each subdomain (see \cref{fig:exam1}),
leading to the numerical second-order convergence rates presented in  \cref{Intersect:table1}. In the other 6 examples, we increase the frequency of the coefficient $a$. The discontinuities at intersection points that are far from the boundary may result in the solution 
$u$ losing regularity such that the numerical convergence rates of FEM or FDM
are relatively small in \cref{Intersect:table2,Intersect:table3,Intersect:table4,Intersect:table5,Intersect:table3:New,Intersect:table7}.

From the performance of FEM and FDM in  \cref{Intersect:table3,Intersect:table4,Intersect:table5} of \cref{Intersect:ex3,Intersect:ex4,Intersect:ex5}, we can observe that the convergence rates of FEM are greater than those of 
FDM. The reason may be that the coefficient functions contain both high frequencies and large jumps in \cref{Intersect:ex3,Intersect:ex4,Intersect:ex5} such that  $u\in H^{1+\lambda}(P_\Omega)$ with a small $\lambda$. 
Another possible reason is that the FDM scheme does not include the intersection grid point (see \cref{fig:interface,ignore:intersection}), i.e., the $(u^{D}_h)_{i,j}$ at the intersection grid point $(x_i,y_j)$ is ignored.
\section{Conclusion}\label{sec:Conclu}
In this paper, we apply second-order finite element method (FEM) and finite difference method (FDM) to numerically solve elliptic cross-interface problems \eqref{intersect:3} with piecewise constant coefficients, homogeneous jump conditions, continuous source terms, and Dirichlet boundary conditions. The numerical results reveal interesting phenomena:
\begin{itemize}
	\item
When the coefficient functions only exhibit minor jumps or low-frequency oscillations, the numerical solutions obtained from FEM and FDM are similar (see \cref{fig:exam1,fig:exam2,fig:exam3:New}).

	\item In contrast, when solving \eqref{intersect:3} with high-contrast and high-frequency coefficient functions,  FEM and FDM produce markedly different solutions  (see \cref{fig:exam3,fig:exam4,fig:exam5,fig:exam7}).
\end{itemize}
To the best of our knowledge, there is currently no literature that has revealed such distinct differences in the numerical solutions obtained from finite element and finite difference methods.
These findings are particularly relevant in the context of the SPE10 benchmark problem, which contains high-contrast and high-frequency permeability resulting from the complicated geological structure in porous media. This discrepancy between FEM and FDM solutions should be carefully considered, especially when developing multiscale  methods where reference solutions are often obtained by the standard FEM. We provide all necessary details to facilitate the replication of our numerical results which allows readers to easily verify the accuracy of our observations.

We only find these interesting numerical results in this paper and do not theoretically determine the reason or provide the rigorous analysis. We can only conclude that the standard FEM and FDM are insufficient to obtain the solution when the coefficient function $a$
includes high jumps and high frequencies. Therefore, it is necessary to analyze the regularity of the solution in the model problem \eqref{intersect:3}, which features high-contrast and high-frequency coefficient functions, and to design reliable FEM and FDM schemes to provide convincing solutions.


\begin{thebibliography}{99}
	
	
	
	
	\bibitem{AliMankad2020} A.~Ali, H.~Mankad, F.~Pereira, and F.~S.~Sousa,  The multiscale perturbation method for second order elliptic equations.
	\emph{Appl. Math. Comput.} \textbf{387} (2020), 125023.
	
	
	

	
	\bibitem{ArbogastTao2013} T.~Arbogast, Z.~Tao, and H.~Xiao, Multiscale mortar mixed methods for heterogeneous elliptic problems.
	\emph{Contemp. Math.} \textbf{586} (2013), 9-21.
	
	
	
	
	\bibitem{ArbogastXiao2013} T.~Arbogast and H.~Xiao,  A Multiscale Mortar Mixed Space Based on Homogenization for Heterogeneous Elliptic Problems.
	\emph{SIAM J. Numer. Anal.} \textbf{51} (2013), 377-399.
	
	
	
	
	\bibitem{TArbHXiao2015} T.~Arbogast and H.~Xiao,  Two-level mortar domain decomposition preconditioners for heterogeneous elliptic problems.
	\emph{Comput. Methods Appl. Mech. Engrg.} \textbf{292} (2015), 221-242.
	
	
	
	\bibitem{Blumen85}
	M.~Blumenfeld,
	The regularity of interface-problems on corner-regions.	\emph{Lecture Notes in Math.} \textbf{1121} (1985), 38-54.
	
	
	
	
	\bibitem{Butler2020} R.~Butler, T.~Dodwell, A.~Reinarz, A.~Sandhu, R.~Scheichl, and L.~Seelinger,  High-performance dune modules for solving large-scale, strongly anisotropic elliptic problems with applications to aerospace composites.
	\emph{Comput. Phys. Commun.} \textbf{249} (2020), 106997.
	
	
	
	
	
	
	\bibitem{CaKi01} Z.~Cai and S.~Kim,
	A finite element method using singular functions for the Poisson equation: corner singularities.
	\emph{SIAM J. Numer. Anal.} \textbf{39} (2001), no. 1,  286-299.
	
	

				
	\bibitem{ChenDai2002}
	Z.~Chen and S.~Dai
	On the efficiency of adaptive finite element methods for elliptic problems with discontinuous coefficients. \emph{SIAM J. Sci. Comput.}     \textbf{24}  (2002),  443-462.
	
	
	
	\bibitem{CFL19}
	X.~Chen, X.~Feng, and Z.~Li,
	A direct method for accurate solution and gradient
	computations for elliptic interface problems. \emph{Numer. Algorithms.}     \textbf{80} (2019), 709-740.
	
	

	
	\bibitem{DFL20}
	B.~Dong, X.~Feng, and Z.~Li,
	An FE-FD method for anisotropic elliptic interface problems. \emph{SIAM J. Sci. Comput.}     \textbf{42}  (2020),  B1041-B1066.
	
	

	
	
	
	
	\bibitem{EngwerHenning2019} C.~Engwer, P.~Henning, A.~M\aa lqvist, and D.~Peterseim,
	Efficient implementation of the localized orthogonal decomposition method. \emph{Comput. Methods Appl. Mech. Engrg.} \textbf{350} (2019), 123-153.
	
	
	
	\bibitem{EwingLLL99}
	R.~Ewing, Z.~Li, T.~Lin, and Y.~Lin,
	The immersed finite volume element methods for the elliptic interface problems. \emph{Math. Comput. Simul.} \textbf{50} (1999), 63-76.
	
	

	
	
\bibitem{FHM21Helmholtz} Q.W.~Feng, B.~Han, and M.~Michelle, 	Sixth-order compact finite difference method for 2D Helmholtz equations with singular sources and reduced pollution effect.    \emph{Commun.Comput. Phys.}  \textbf{34} (2023),  672-712. 



\bibitem{FengHanMinevPOISS}
Q.W.~Feng, B.~Han, and P.~Minev, Sixth order compact finite difference schemes for Poisson interface problems with singular sources. \emph{Comput. Math. Appl.} \textbf{99} (2021), 2-25.




\bibitem{FengHanMinevCross} Q.W.~Feng, B.~Han, and P.~Minev,
Compact  9-point  finite difference methods with high accuracy order and/or  M-Matrix property
for  elliptic cross-interface problems.  \emph{J. Comput. Appl. Math.}  \textbf{428}  (2023), 115151.




\bibitem{FengHanMinevHYB} Q.W.~Feng, B.~Han, and P.~Minev,
Sixth-order hybrid  finite difference methods for elliptic interface problems with mixed boundary conditions.   \emph{J. Comput. Phys.} 497 (2024) 112635. 




	\bibitem{FendZhao20pp109677}
	H.~Feng and S.~Zhao, A fourth order finite difference method for solving elliptic
	interface problems with the FFT acceleration. \emph{J. Comput. Phys.}     \textbf{419} (2020), 109677.	
	
	
	
     \bibitem{FGW73} G.~J.~Fix, S.~Gulati and G.~I.~Wakoff,
	On the use of singular functions with finite element approximations.
	\emph{J. Comput. Phys.} \textbf{13} (1973), 209-228.
	
	
	
	
	\bibitem{FuChungLi2019} S.~Fu, E.~Chung, and G.~Li, Edge multiscale methods for elliptic problems with heterogeneous coefficients. \emph{J. Comput. Phys.} \textbf{396} (2019),  228-242.
	
	
	
	
	\bibitem{GongLiLi08} Y.~Gong, B.~Li, and Z.~Li,
	Immersed-interface finite-element methods for elliptic interface problems with nonhomogeneous jump conditions.
	\emph{SIAM J. Numer. Anal.} \textbf{46} (2008), 472-495.	
	
	
	
	\bibitem{GuiraAusas2019} R.~T.~Guiraldello, R.~F.~Ausas, F.~S.~Sousa, F.~Pereira, and G.~C.~Buscaglia, Interface spaces for the Multiscale Robin Coupled Method in reservoir simulation.
	\emph{Math. Comput. Simul.} \textbf{164} (2019), 103-119.
	
	
	

	
	\bibitem{HeLL2011} X.~He, T.~Lin, and Y.~Lin,
	Immersed finite element methods for elliptic interface problems with non-homogeneous jump conditions.
	\emph{Int. J. Numer. Anal. Model.} \textbf{8} (2011), 284-301.
	
	
	\bibitem{Hellmanvist2017} F.~Hellman and A.~M\aa lqvist,
Contrast independent localization of multiscale problems. \emph{SIAM Multiscale Model. Simul.} \textbf{15} (2017), 1325-1355.

	
	

\bibitem{HouHwang2017} T.~Y.~Hou, F.~N.~Hwang, P.~Liu, and C.~C.~Yao, An iteratively adaptive multi-scale finite element method for elliptic PDEs with rough coefficients. \emph{J. Comput. Phys.} \textbf{336} (2017),  375-400.



	
	\bibitem{Jaramillo2022} A.~Jaramillo, R.~T.~Guiraldello, S.~Paz, R.~F.~Ausas, F.~S.~Sousa, F.~Pereira, and G.~C.~Buscaglia, Towards HPC simulations of billion-cell reservoirs by multiscale mixed methods.
	\emph{Comput. Geosci.} \textbf{26} (2022), 481-501.
	
	
	
	
		\bibitem{Kell71} R.~B.~Kellogg,
	Singularities in interface problems.
	\emph{Numerical Solution of Partial Differential Equations-\uppercase\expandafter{\romannumeral2}} (1971), 351-400.
	
	
	\bibitem{Kell72} R.~B.~Kellogg,
	Higher order singularities for interface problems.
	\emph{The Mathematical Foundations of the Finite Element Method with Applications to Partial Differential Equations} (1972), 589-602.
	
	
	
	\bibitem{Kell75}  R.~B.~Kellogg,
	On the Poisson equation with intersecting interfaces.
	\emph{Appl. Anal.} \textbf{4} (1975), no. 2,  101-129.
	
	
	
		\bibitem{KCPK07} S.~Kim, Z.~Cai, J.~Pyo and S.~Kong,
	A finite element method
	using singular functions: interface problems.
	\emph{Hokkaido Math. J.} \textbf{36} (2007), no. 4, 815-836.
	
	
	\bibitem{KiKo06} S.~Kim and S.~Kong,
	A finite element method dealing the singular points with a cut-off function.
	\emph{J. Appl. Math. Comput.} \textbf{21} (2006), no. 1-2,  141-152.
	
	
	
	
	\bibitem{Kippe2008} V.~Kippe,  J.~E.~Aarnes, and K.~A.~Lie, A comparison of multiscale methods for elliptic problems in porous media flow.
	\emph{Comput. Geosci.} \textbf{12} (2008), 377-398.
	
	
	
	
	
	
	\bibitem{LeLi94}
	R.~J.~ Leveque and Z.~Li,
	The Immersed interface method for elliptic equations with discontinuous coefficients and singular sources. \emph{SIAM J. Numer. Anal.}     \textbf{31} (1994),  1019-1044.
	
	
	
	
	
	\bibitem{Li98}
	Z.~Li,
	A fast iterative algorithm for elliptic interface problems. \emph{SIAM J. Numer. Anal.}     \textbf{35} (1998),  230-254.
	
	
	
	
	
	\bibitem{LiHu2021} G.~Li and J.~Hu, Wavelet-based edge multiscale parareal algorithm for parabolic equations with heterogeneous coefficients and rough initial data. \emph{J. Comput. Phys.} \textbf{444} (2021),  110572.
	
	
	
	
		
	
	\bibitem{lqvistPersson2018} A.~M\aa lqvist and A.~Persson,
	Multiscale techniques for parabolic equations. \emph{Numer. Math.} \textbf{138} (2018), 191-217.
	
	
	
	
	
			
	\bibitem{lqvistPeterseim2014} A.~M\aa lqvist and D.~Peterseim,
	Localization of elliptic multiscale problems. \emph{Math. Comput.} \textbf{83} (2014), 2583-2603.
	
	
	

	
	
	
	\bibitem{Minev2018}
	P.~Minev, S.~Srinivasan, and P.~N.~Vabishchevich, Flux formulation of parabolic equations with highly heterogeneous coefficients. \emph{J. Comput. Appl. Math.} \textbf{340} (2018), 582-601.
	
	
	
	
	\bibitem{Nic90} S.~Nicaise,
	Polygonal interface problems: higher regularity
	results.
	\emph{Commun. Partial. Differ. Equ.} \textbf{15} (1990), no. 10, 1475-1508.
	
	
	\bibitem{NS94} S.~Nicaise and  A.~M.~Sandig,
	General interface problems-\uppercase\expandafter{\romannumeral1}.
	\emph{Math. Methods Appl. Sci.} \textbf{17} (1994), 395-429.
	
	

	
	
	\bibitem{PanHeLi21} K.~Pan, D.~He, and Z.~Li,
	A high order compact FD framework for elliptic BVPs
	involving singular sources, interfaces, and irregular domains,
	\emph{J. Sci. Comput.} \textbf{88} (2021), 1-25.
	
	
	
	
	
	\bibitem{Petz01} M.~Petzoldt,
	Regularity results for Laplace interface problems in two dimensions.
	\emph{J. Anal. Appl.} \textbf{20} (2001), no. 2,  431-455.
	
	\bibitem{Petz02} M.~Petzoldt,
	A posteriori error estimators for elliptic equations with
	discontinuous coefficients.
	\emph{Adv. Comput. Math.} \textbf{16} (2002), no. 1, 47-75.
	
	
	
	
	\bibitem{PietroErn2012} D.~A.~Di Pietro and A.~Ern, Mathematical aspects of discontinuous Galerkin methods. \emph{Springer Berlin, Heidelberg.} 2012.
	
	
	
	\bibitem{Rasaei2008} M.~R.~Rasaei and M.~Sahimi, Upscaling and simulation of waterflooding in heterogeneous reservoirs using wavelet transformations: application to the SPE-10 model.
	\emph{Transp. Porous. Med.} \textbf{72} (2008), 311-338.
	
	
	
	
	

	
	
	\bibitem{Tahmasebi2018} P.~Tahmasebi and S.~Kamrava, A multiscale approach for geologically and flow consistent modeling.
	\emph{Transp. Porous. Med.} \textbf{124} (2018), 237-261.
	
	
	
	
	
	\bibitem{Vazqu07} J.~L.~V{\'a}zquez, The Porous medium equation: mathematical theory. \emph{Clarendon Press.} 2007.
	
	
	
	
	\bibitem{WieBube00} A.~Wiegmann and  K.~P.~Bube,
	The explicit-jump immersed interface method: finite difference methods for PDEs with piecewise smooth solutions. \emph{SIAM J. Numer. Anal.}     \textbf{37} (2000),    827-862.
	
	
	
	
	\bibitem{YuWei073D} S.~Yu and G.~W.~Wei,
	Three-dimensional matched interface and boundary
	(MIB) method for treating geometric singularities.
	\emph{J. Comput. Phys.} \textbf{227} (2007), 602-632.
	
	
	\bibitem{YuZhouWei07} S.~Yu, Y.~Zhou,  and G.~W.~Wei,
	Matched interface and boundary (MIB) method for
	elliptic problems with sharp-edged interfaces. \emph{J. Comput. Phys.}     \textbf{224} (2007),  729-756.
	
	
	
		
	\bibitem{ZhangCiHou2015} Z.~Zhang, M.~Ci, and T.~Y.~Hou,
	A multiscale data-driven stochastic method for elliptic PDEs with random coefficients. \emph{SIAM Multiscale Model. Simul.} \textbf{13} (2015), 173-204.
	
	
	
	\bibitem{XiaolinZhong07} X.~Zhong,
	A new high-order immersed interface method for solving
	elliptic equations with imbedded interface of discontinuity. \emph{J. Comput. Phys.} \textbf{225} (2007),  1066-1099.
	
	
	
	\bibitem{ZW06} Y.~C.~Zhou and G.~W.~Wei,
	On the fictitious-domain and interpolation formulations of
	the matched interface and boundary (MIB) method. \emph{J. Comput. Phys.} \textbf{219} (2006),  228-246.
	
	
	
	\bibitem{ZZFW06} Y.~C.~Zhou, S.~Zhao, M.~Feig, and G.~W.~Wei,
	High order matched interface and boundary method for elliptic
	equations with discontinuous coefficients and singular sources. \emph{J. Comput. Phys.} \textbf{213} (2006), 1-30.
	


	










	






\end{thebibliography}
\end{document}